\documentclass[12pt]{oxarticle}
\pdfoutput=1
\usepackage{graphicx, float, array, xspace, amscd, amsmath, amsthm, amssymb}
\usepackage[sort&compress, comma, square, numbers]{natbib}
\usepackage{cjhebrew}

\let\SS=\S 
\let\textbreve=\u
\renewcommand{\a}{\alpha}
\renewcommand{\b}{\beta}
\newcommand{\g}{\gamma}
\renewcommand{\d}{\delta}
\newcommand{\e}{\epsilon}\newcommand{\ve}{\varepsilon}
\newcommand{\z}{\zeta}

\renewcommand{\k}{\kappa}

\newcommand{\x}{\xi}

\newcommand{\p}{\pi}
\renewcommand{\r}{\rho}
\newcommand{\s}{\sigma}\renewcommand{\S}{\Sigma}
\renewcommand{\t}{\tau}
\renewcommand{\u}{\upsilon}\newcommand{\U}{\Upsilon}
\newcommand{\vph}{\varphi}

\newcommand{\ps}{\psi}
\renewcommand{\o}{\omega}\renewcommand{\O}{\Omega}

\newcommand{\cA}{\mathcal{A}}
\newcommand{\cB}{\mathcal{B}}

\newcommand{\cG}{\mathcal{G}}

\newcommand{\cK}{\mathcal{K}}

\newcommand{\cM}{\mathcal{M}}

\newcommand{\cO}{\mathcal{O}}
\newcommand{\cP}{\mathcal{P}}

\newcommand{\cS}{\mathcal{S}}

\newcommand{\cX}{\mathcal{X}}

\newcommand{\IC}{\mathbb{C}}

\newcommand{\IF}{\mathbb{F}}

\newcommand{\IP}{\mathbb{P}}

\newcommand{\IR}{\mathbb{R}}

\newcommand{\IZ}{\mathbb{Z}}

\font\eightrm=cmr8 at 8pt
\font\csc=cmcsc10

\newcommand{\beq}{\begin{equation}}
\newcommand{\eeq}{\end{equation}}
\newcommand{\beqnn}{\begin{equation*}}
\newcommand{\eeqnn}{\end{equation*}}
\newcommand{\bea}{\begin{eqnarray}}
\newcommand{\eea}{\end{eqnarray}}
\newcommand{\bean}{\begin{eqnarray*}}
\newcommand{\eean}{\end{eqnarray*}}

\newcommand{\fref}[1]{Figure~\ref{#1}}
\newcommand{\tref}[1]{Table~\ref{#1}}
\newcommand{\sref}[1]{\SS\ref{#1}}

\newcommand{\defineas}{\buildrel\rm def\over =}

\newcommand{\ord}[1]{\cO\kern-2pt\left(#1\right)}

\newcommand{\Fp}{\IF_p}

\newcommand{\ii}{\text{i}}

\newcommand{\place}[3]{\vbox to0pt{\kern-\parskip\kern-7pt
                             \kern-#2truein\hbox{\kern#1truein #3}
                             \vss}\nointerlineskip}

\newcommand{\smallfrac}[2]{\frac{\scriptstyle #1}{\scriptstyle #2}}

\DeclareFontFamily{U}{wncy}{}
\DeclareFontShape{U}{wncy}{m}{n}{<->wncyr10}{}
\DeclareSymbolFont{mcy}{U}{wncy}{m}{n}
\DeclareMathSymbol{\sha}{\mathord}{mcy}{"58}

\newcommand{\capt}[3]{\parbox{#1}{\renewcommand{\baselinestretch}{1.0}
                                                           \caption{\label{#2}\small\it #3}}}

\newcommand{\dP}{\hbox{dP}}
\newcommand{\tC}{\widetilde{C}}
\newcommand{\Aut}{\text{Aut}}
\newcommand{\Pic}{\text{Pic}}

\newcommand{\Mustata}{Musta\c{t}\textbreve{a}\xspace}

\renewcommand{\baselinestretch}{1.1}
\numberwithin{equation}{section}
\proofmodefalse
\begin{document}
\pagestyle{empty}

\begin{center}
\null\vskip0.3in
{\Huge Lines on the Dwork Pencil\\
of\\
Quintic Threefolds\\[0.4in]}
{\csc Philip Candelas$^1$, Xenia de la Ossa$^1$\\[0.1in]
Bert van Geemen$^2$ and Duco van Straten$^3$\\[0.4in]}
{\it $^1$Mathematical Institute\hphantom{$^1$}\\
University of Oxford\\
24-29 St.\ Giles'\\
Oxford OX1 3LB, UK\\[0.3in]}
\parbox{5.6in}{\it
\begin{minipage}{2.5in}
\centering
$^2$Dipartimento di Matematica\hphantom{$^1$}\\ 
Universit\`{a} di Milano\\ 
Via Cesare Saldini 50\\
20133 Milano, Italy\\
\end{minipage}
\hskip0.5in
\begin{minipage}{2.5in}
\centering
$^3$Fachbereich 17\hphantom{$^1$}\\ 
AG Algebraische Geometrie\\
Johannes Gutenberg-Universit\"at\\
D-55099 Mainz, Germany\\
\end{minipage}
}
\vfill
{\bf Abstract}\\[10pt]
\begin{minipage}{\textwidth}
\small
We present an explicit parametrization of the families of lines of the Dwork pencil of quintic threefolds. This gives rise to isomorphic curves $\tC_{\pm\vph}$ which parametrize the lines. These curves are 125:1 covers of genus six curves $C_{\pm\vph}$. The $C_{\pm\vph}$ are first presented as curves in $\IP^1{\times}\IP^1$ that have three nodes. It is natural to blow up $\IP^1{\times}\IP^1$ in the three points corresponding to the nodes in order to produce smooth curves. The result of blowing up $\IP^1{\times}\IP^1$ in three points is the  quintic del Pezzo surface $\dP_5$, whose automorphism group is the permutation group $\cS_5$, which is also a symmetry of the pair of curves $C_{\pm\vph}$. The subgroup $\cA_5$, of even permutations, is an automorphism of each curve, while the odd permutations interchange $C_\vph$ with $C_{-\vph}$. The ten exceptional curves of $\dP_5$ each intersect the $C_\vph$ in two points corresponding to van Geemen lines. We find, in this way, what should have anticipated from the outset, that the curves $C_\vph$ are the curves of the Wiman pencil.  We consider the family of lines also for the cases that the manifolds of the Dwork pencil become singular. For the conifold the curve $C_\vph$ develops six nodes and may be resolved to a~$\IP^1$. The group $\cA_5$ acts on this $\IP^1$ and we describe this action.
\end{minipage}

\end{center}
\newpage
{\renewcommand{\baselinestretch}{1.0}\tableofcontents}
\newpage
\setcounter{page}{1}
\pagestyle{plain}

\begin{minipage}{6.5in}
\begin{cjhebrew}\raggedleft\large
ye+s d*AbAr +sEy*o'mar r:'eh zEh .hAdA+s hU'\\ 
k*:bAr hAyAh l:`olAmiym 'a:+sEr hAyAh mil*:pAnenU;\\[10pt]
\end{cjhebrew}
\raggedleft\it
Is there any thing whereof it may be said, See, this is new?\\
it hath been already of old time, which was before us.\\[5pt]
Ecclesiastes 1:10
\end{minipage}
\vspace{20pt}

\section{Introduction}
\subsection{Lines on the cubic surface and quintic threefold}
This article concerns the lines contained in the Dwork pencil of quintic threefolds. These manifolds, which we denote by $\cM_\ps$, are realised as hypersurfaces in $\IP^4$ by the quintics
\beq
\sum_{j=1}^5 x_j^5 - 5 \psi\, x_1x_2x_3x_4x_5 ~=~0~.
\label{DworkPencil}\eeq
The study of the lines on quintic threefolds has a history going back to Schubert in the 19th century, who calculated that the {\em generic} quintic contains 2875 lines, in fact Schubert performed the calculation twice, using different methods \cite{Schubert1, Schubert2}. The quintics of the Dwork pencil are, however, far from being generic and are known to contain continuous families of lines. 

Before summarizing the history of our understanding of lines on the quintic it might be useful to recall that this study began as a natural extension of the classical study of lines on cubic surfaces. These lines were discovered by Cayley and Salmon. The story is famous: Cayley remarked in a letter that counting constants suggested a finite number, and Salmon gave immediately the number 27 in response to the letter. The results of this correspondence were published in 1849~\cite{CayleyLines, SalmonLines}.
The configuration of the lines and their intricate symmetries have been of topic of fascination to algebraic geometers ever since. A classical source of information is the book of Henderson~\cite{HendersonLines}.

There are differences between the cubic and the quintic; in order to appreciate these let us recall the most elementary facts. The Fermat cubic in $\IP^3$ is given by the equation
\beq
\sum_{r=1}^4 y_r^3~=~0~.
\label{FermatCubic}\eeq
This surface contains the lines $y_r = (u, -\o^j u, v, -\o^k v)$, where $\o$ denotes a nontrivial cube root of unity and $1\leq j,k\leq3$. By permuting the coordinates we find 27 lines that lie in the cubic and this is the total number. The beautiful and surprising fact is that if we deform the cubic, the lines deform with the surface so that there are always 27 lines. For a generic cubic it will be hard to see the lines explicitly. In fact 
C. Jordan~\cite{Jordan1870} showed that the Galois group on which the determination of the lines depends is in general a simple group of order 25,920 which can be identified with the Weyl group of the lattice $E_6$ modulo its center. A modern reference for these results is~\cite{HarrisGaloisGroups}.

For the Dwork pencil \eqref{DworkPencil} the situation is already more complicated, even for the case of the Fermat quintic with $\psi=0$. For this case we may write down analogues of the lines that exist for the Fermat cubic
$$
x_j~=~(u, -\z^k u, v, -\z^\ell v,0)~,
$$
with $\z$ a nontrivial fifth root of unity and $1\leq k,\ell\leq5$. By permuting coordinates and taking all values of $k$ and $\ell$ we find 375 such lines, which we will refer to here as the isolated lines\footnote{These lines are often known as the exceptional lines, however, to refer to them as such here would invite confusion with the exceptional lines of the del Pezzo surface $\dP_5$, to which we shall make frequent reference. These lines are indeed isolated for $\psi\neq 0$, but, as we shall see, they lie in continuous families of lines for~$\psi=0$.}. 
Note that, since one of the coordinates vanishes identically, these lines lie in $\cM_\psi$ for all $\psi$. 

There are other lines also. Consider those of the form
$$
x_j~=~(u,-\z^k u, av, bv, cv)~~~\text{with}~~~a^5 + b^5 + c^5 ~=~ 0~.
$$
For given $k$, these give rise to a cone of lines, that all pass through the point $(1, -\z^k,0,0,0)$, and are parametrized by the curve $a^5+b^5+c^5=0$ in $\IP^2$. By counting the different values of~$k$ and the inequivalent permutations of the coordinates we see that there are 50 cones of lines. The cones contain the isolated lines. In fact, the isolated lines are the lines in which the cones meet. For example
the cones $(u,-u,av,bv,cv)$ and $(\tilde{a}u,\tilde{b}u,v,-v,\tilde{c}u)$ meet in the isolated line 
$(u,-u,v,-v,0)$. Each cone contains 15 isolated lines and meets 15 other cones in these lines. If two cones intersect, they do so in precisely one of the isolated lines.

In~\cite{MR1024767} it is shown that there are no further lines in $\cM_0$ beyond the cones and the isolated lines and, furthermore, that, under a sufficiently general deformation, each isolated line splits into 5 lines and each cone breaks up into 20 discrete lines, yielding the correct total of 
$50{\times}20 + 5{\times}375 = 2875$ discrete lines. 

A quintic threefold deforms with 101 parameters, and for generic values of these parameters there are, as has been observed, 2875 discrete lines. It is known, however, that there are families of quintic threefolds that deform with 100 parameters, for which the configuration of lines is degenerate~\cite{KatzDegenerations}.

Let us return now to the one parameter family $\cM_\psi$ for $\psi\neq0$. The manifolds of the Dwork pencil have a large group of automorphisms isomorphic to $\cS_5{\rtimes}\cG$, where $\cS_5$ is the permutation group acting on the five coordinates and $\cG\cong (\IZ/5\IZ)^3$ has the action
$$
(x_1,\, x_2,\, x_3,\, x_4,\, x_5)\longrightarrow (\z^{n_1}\, x_1,\, \z^{n_2}\, x_2,\, \z^{n_3}\, x_3,\,
\z^{n_4}\,  x_4,\, \z^{n_5}\, x_5)~~~\text{with}~~~\sum_{j=1}^5 n_j=0\bmod 5~.
$$

In the 1980's one of the present authors (BvG) found special lines that lie in the $\cM_\psi$. These eponymous lines are important in what follows so we shall pause, presently, to review their properties. For the moment we simply note that there are 5,000 such lines, so since this number exceeds 2875, there must be, possibly in addition to discrete lines, a continuous family~\cite{MR1085631}. It was subsequently proved by 
Anca \Mustata~\cite{Mustata:fk}, using sophisticated methods, that, for $\ps\neq 0$, $\cM_\psi$ contains two continuous families of lines, parametrized by isomorphic curves, $\tC_\pm$, of genus 626, and the 375 isolated lines as the only lines that do not lie in the continuous families. The genus 626 curves have Euler number $\chi=2-2{\times}626=-1250$. It follows from the theory of the Abel-Jacobi mapping (see some further remarks in \sref{ABcurves}) that under a generic deformation, each of these curves gives rise to 1250 discrete lines, so that, all together, there are again $375 + 2{\times}1250=2875$ lines.

One of our aims here is to parametrize the two families of lines, $\tC_\pm$, explicitly. The surprise is that the explicit parametrization is not as complicated as might have been anticipated. 
\subsection{The van Geemen lines}
If the $\cM_\ps$ were to contain 2875 lines `as expected' we would want to find the $2875\, {-}\, 375 = 2500$ lines that are missing (assuming that the special lines are to be counted with multiplicity one). Now $\cS_5$ has subgroups of order three, for example the subgroup that permutes $(x_2,x_4,x_5)$ cyclically (the reason for choosing this particular subgroup is to conform with a choice of parametrisation that will come later). The number of missing lines is not divisible by three so some would have to be fixed (as lines but not necessarily pointwise) by the subgroup. This motivates seeking lines that are invariant under the proposed subgroup.

The points that are invariant under the subgroup are of the form
$$
(a,d,b,d,d)~,~~~(0,1,0,\o,\o^2)~, ~~~(0,1,0,\o^2,\o)~.
$$
It is immediate that the plane $(a,d,b,d,d)$ does not contain a line of $\cM_\ps$ and that the line passing through 
$(0,1,0,\o,\o^2)$ and $(0,1,0,\o^2,\o)$ does not lie in $\cM_\ps$. Consider however the line
\beq
u\,\big(1,d,b,d,d\big) + (v - d u)\,\big(0,1,0,\o,\o^2 \big)~=~\big(u,v,bu,cu+\o v, -\o^2(cu-v)\big)~,
\label{VanGzero}\eeq
where $c=(1-\o)d$. This line lies in $\cM_\ps$ provided
\beq 
b~=~\frac32\, \ps\g^2~,~~~c~=~\frac12\,(1-\o)\ps\,\g~, 
\label{VanGCondone}\eeq
with $\g$ a solution of the tenth order equation
\beq
\g^{10}-\frac19\,\g^5+\left(\frac{2}{3\ps}\right)^5=0~.
\label{VanGCondtwo}\eeq

Given that the lines \eqref{VanGzero}, subject to \eqref{VanGCondone} and \eqref{VanGCondtwo} lie in 
$\cM_\psi$ it is clear that so do lines of the form
\beq
\Big(u,\, v,\, \z^{-k-\ell}\, b u,\, \z^k (cu+\o v),\, -\z^\ell \o^2 (cu - v)\Big) ~,
\label{VanGone}\eeq
with $\z$ is a nontrivial fifth root of unity, $1\leq k,\ell\leq 5$, since these are images of the previous line under the action of $\cG$. The van Geemen lines are the lines that are equivalent to this more general form, up to permutation of coordinates. These, more general, lines are no longer invariant under the cyclic permutation of 
$(x_2,x_4,x_5)$. However, since they are in the $\cS_5{\rtimes}\cG$ orbit of \eqref{VanGzero}, which has an 
$\cS_5$ stabilizer of order three, the more general lines each have a stabilizer of order three.

There are changes of coordinates that preserve the general form of a van~Geemen line. Setting
$u=\z^{k+\ell}\tilde{u}/b$ effectively interchanges the $u$ and $bu$ terms by bringing the line                    \eqref{VanGone} to the form 
$$
\Big(\, \z^{-k-\ell}\, \tilde{b}\tilde{u},\, v,\, \tilde{u},\, \z^k (\tilde{c}\tilde{u}+\o v),\, 
-\z^\ell \o^2 (\tilde{c}\tilde{u} - v)\,\Big) 
$$
where
$$
\tilde{b}=\frac{\z^{2(k+\ell)}}{b}=\frac32\ps\,\tilde{\g}^2~~~\hbox{and}~~~
\tilde{c}=\z^{k+\ell}\,\frac{c}{b} = \frac12(1-\o)\ps\,\tilde{\g}~~~\hbox{with}~~~
\tilde{\g}=\z^{k+\ell}\frac{2}{3\ps\g}~,
$$
and in these relations $\tilde{\g}$ is another root of equation \eqref{VanGCondtwo}.

If we return to \eqref{VanGone} and write 
$$
v_1=v~,~~~v_2= cu+\o v~,~~~v_3=-\o^2(cu-\o v) 
$$
and change coordinates and parameters by setting
$$
\tilde{v} = \z^k v_2~,~~~\tilde{b}=\z^{2k}b~,~~~\tilde{c}=\z^k c
$$
then we have
$$
\tilde{v}_1\defineas \tilde{v}=\z^k v_2~,~~~\tilde{v}_2\defineas \tilde{c}u+\o\tilde{v}=\z^k v_3~,~~~
\tilde{v}_3\defineas -\o^2(\tilde{c}u-\tilde{v})=\z^k v_1
$$
and the effect of the coordinate transformation is
$$
(u,\, v_1,\, \z^{-k-\ell}bu,\, \z^k v_2,\, \z^\ell v_3) = 
(u,\, \z^{-k}\tilde{v}_3,\,  \z^{2k-\ell}\tilde{b}u,\, \tilde{v}_1,\, \z^{\ell-k}\tilde{v}_2)~.
$$
Note that the change in $b$ and $c$ is consistent with $\g\to\tilde{\g}=\z^k\g$ and $\tilde{\g}$ is another root of~\eqref{VanGCondtwo}. In this way one may, in effect, rotate the quantities $v_j$ cyclically, however we are left with two orderings of the $v_j$ that cannot be transformed into each other. 

The counting is that, up to coordinate redefinitions, there are 10 ways to choose two positions for the components $u$ and $bu$ and a further two choices in the placing of the components $v_j$. There are two choices for $\o$, five for $\g$, given $\g^5$, and 25 ways to choose $k$ and $\ell$. Thus there are, in total,
$10{\times}2{\times}2{\times}5{\times}25=5,000$ van Geemen lines. In this accounting we consider 
\eqref{VanGCondtwo} to be a quadratic equation for $\g^5$ and we do not count the two roots separately since these are interchanged by the coordinate transformation that interchanges $u$ and $bu$. The fact that there are 5,000 van Geemen lines while $\#(\cS_5{\rtimes}\cG)=5!{\times}5^3=15,000$ again implies (though one can also check this directly) that each of these lines has a stabilizer of order exactly~three. 

Since the number of lines, if discrete, must be 2875, counted with multiplicity, the fact that 5000 lines have been identified implies that, while there may be discrete lines,  there must also be a continuous family of lines.

If we pick a particular value for $\g$ and act with an element of $\cG$
as above on the line
$$
\big(u,\, v,\, b u,\, cu+\o v,\, - \o^2 (cu - v)\big)$$
and then set $u=\z^{-n_1}\tilde{u}$,  $v=\z^{-n_2}\tilde{v}$,  $\g=\z^{n_1-n_2}\tilde{\g}$ and make the corresponding changes $b=\z^{2(n_1-n_2)}\tilde{b}$ and $c=\z^{n_1-n_2}\tilde{c}$ then we obtain the line
$$
(\tilde{u},\, \tilde{v},\, \z^{n_1-2n_2+n_3}\tilde{b}\tilde{u},\,
\z^{n_4-n_2}(\tilde{c}\tilde{u}+\o\tilde{v}),\, -\z^{n_5-n_2}\o^2(\tilde{c}\tilde{u}-\tilde{v}))~.
$$
In this way we obtain 125 copies of a van Geemen line by acting with $\cG$ on a particular line, provided that we understand $\cG$ to act on $\g$ as indicated.
\subsection{The Wiman pencil}
In 1897 Wiman~\cite{Wiman} noted the existence of a remarkable plane sextic curve $C_0$, with four nodes, that is invariant under the permutation group $\cS_5$. These automorphisms appeared the more mysterious owing to the fact that, of the 120 automorphisms, 96 are realised nonlinearly. The story was taken up by 
Edge~\cite{Edge} after some eighty years, who noted that $C_0$ is ``only one, though admittedly the most interesting'' of a one parameter family of four-nodal sextics $C_\vph$ on which the group $\cS_5$ acts. The action is such that the subgroup $\cA_5$, of even permutations, preserves each $C_\vph$ while the odd permutations interchange $C_\vph$ with $C_{-\vph}$. The curve $C_0$ is known as the Wiman curve and the one parameter family $C_\vph$ is known as the Wiman pencil. Edge notes also that it is natural to blow up the plane in the four nodes of the curves. One obtains, in this way, smooth curves which, in this introduction, we will also denote by $C_\vph$. These smooth curves live in the quintic del Pezzo surface\footnote{There is difference in convention between mathematicians and physicists in writing $\dP_n$. A physicist tends to mean $\IP^2$ blown up in $n$ points, in general position, while a mathematician often means the del Pezzo surface of degree $n$. In the `mathematician's' convention, which we use here, the surface which results from blowing up $\IP^2$ in $n\leq 8$ points, in general position, is $\dP_{9-n}$.} $\dP_5$.  

With our explicit parametrization of the families of lines $\tC_\pm$, and benefit of hindsight, we find what should have been suspected from the outset: the curves $\tC_\pm$ are 125:1 covers of the curves $C_{\pm\vph}$ of the Wiman pencil. Where the parameter $\vph$ is related to the parameter of the quintic by
$$
\vph^2~=~\frac{32}{\ps^5} - \frac34~.
$$
The remarkable action of $\cS_5$ on the curves of the Wiman pencil is seen to correspond to the symmetry of the configuration of the lines of the Dwork quintics.
\subsection{Layout of this paper}
In \sref{families} we present the explicit parametrization of the families of lines. This gives rise to curves 
$C_{\pm\vph}^0$ whose resolutions have  125:1 covers $\tC_\vph$ which parametrize the lines. The curves $C_{\pm\vph}^0$ are first presented as curves in $\IP^1{\times}\IP^1$ that have three nodes. It is noted that the two curves $C_\vph^0$ and $C_{-\vph}^0$ intersect in the three nodes and in 14 other points. Resolution of the nodes replaces each of the nodes by two points which continue to be points of intersection  of the two curves. Thus there are 20 points of intersection and it is noted that each of these correspond to van Geemen lines. It is natural to blow up $\IP^1{\times}\IP^1$ in the three points corresponding to the nodes in order to produce smooth curves $C_{\pm\vph}$. While it is not the case that $\IP^1{\times}\IP^1$ is $\IP^2$ blown up in a point, it is the case that $\IP^1{\times}\IP^1$ blown up in three points is the same as $\IP^2$ blown up in four points, which is the del Pezzo 
surface $\dP_5$. We review the geometry of $\dP_5$ in \SS\ref{dp5}. The first fact to note is that the automorphism group of $\dP_5$ is the permutation group $\cS_5$. There is also an embedding $\dP_5\hookrightarrow \IP^5$ which is useful owing to the fact that the $\cS_5$ transformations become linear, as automorphisms of $\IP^5$, in this presentation of the surface. The surface $\dP_5$ has 10 exceptional curves. These are the blow ups of the four points of $\IP^2$ together with the six lines that pass through the six pairs of points. Three of these exceptional curves resolve the nodes of $C_\vph^0$ and so intersect the resolved curve in two points. These points correspond, as noted previously, to van Geemen lines. The $\cS_5$ automorphisms permute the 10 exceptional curves so we expect that each of the 10 exceptional curves of $\dP_5$ will intersect $C_\vph$ in two points corresponding to van Geemen lines. Checking that this is indeed so is the subject 
of~\SS\ref{secondlook}. In order to properly understand the intersections of the exceptional curves with the $C_\vph$ we consider the Pl\"ucker coordinates of the lines of the quintic and the embedding 
$\dP_5\hookrightarrow \IP^9$. We give also, in this section, a detailed discussion of the 125:1 cover $\tC_\vph \to C_\vph$.
 
In \SS\ref{singularmanifolds} we turn to the form of the curves $C_\vph$ for the cases $\ps^5=0,1,\infty$ that the manifold $\cM_\ps$ either requires special consideration, for the case $\ps=0$, or is singular. For the conifold there are two values $\vph=\pm 5\sqrt{5}/2$ which correspond to $\psi^5=1$. For these, we find that the curve $C_\vph$ develops six nodes and may be resolved to a $\IP^1$. Thus $\tC_\vph$ is the union of 125
$\IP^1$'s. The group $\cA_5$ acts on each of these and we describe this action. 

A number of technical points are relegated to appendices.
\subsection{The zeta function and the $\cA$ and $\cB$ curves}\label{ABcurves}
It is of interest to study the manifolds $\cM_\ps$ of the Dwork pencil over the finite field $\Fp$. The central object of interest, in this situation, is the $\z$-function. For general $\psi$, that is $\psi^5\neq 0,1,\infty$, this takes the form~\cite{Candelas:2004sk} 
$$
\z_\cM(T,\ps)~=~
\frac{R_{\bf  1}(T,\ps)\, R_\cA(p^\r T^\r,\ps)^\frac{20}{\r}\, R_\cB(p^\r T^\r,\ps)^\frac{30}{\r}}{(1-T)(1-pT)(1-p^2T)(1-p^3T)}~.
$$
In this expression the $R$'s are quartic polynomials in their first argument and, here, 
$\r$ \hbox{$(=1,2~\text{or}~4)$} is the least integer such that $p^\r{-}1$ is divisible by 5.
The quartic $R_{\bf  1}$, for example, has the structure
$$ 
R_{\bf  1}(T,\ps)~=~1 + a_{\bf  1}(\ps)\,T + b_{\bf  1}(\ps)\,pT^2 + a_{\bf  1}(\ps)\,p^3T^3 + p^6 \,T^4 
$$
with $a_{\bf  1}$ and $b_{\bf  1}$ integers that vary with $\ps\in\Fp$. The other factors 
$R_\cA$ and $R_\cB$ have a similar structure. The numerator of the $\z$-function corresponds to the 
Frobenius action on $H^3(\cM_\psi)$. It is intriguing that these factors are related to certain genus 4 Riemann curves $\cA$ and $\cB$. What is meant by this is that there is a genus 4 curve $\cA$, that varies with $\psi$, with $\z$-function satisfying
$$
\z_\cA(T,\ps)~=~\frac{R_\cA(T,\ps)^2}{(1-T)(1-pT)}~,
$$
and there is an analogous relation for another curve $\cB$. The intriguing aspect is that the curves $\cA$ and 
$\cB$ are not directly visible in $\cM_\psi$.

The theory of the Abel-Jacobi mapping provides a context of explaining this phenomenon. More precisely, a loop $\g \in H_1(\tC_{\pm\vph})$ determines a 3-cycle $T(\g) \in \cM_\psi$ which is the union of the lines corresponding to the points of $\g$. By duality one obtains a map 
$a: H^3(\cM_\psi) \to H^1(\tC_{\pm\vph})$, whose
kernel should have dimension 4 and giving rise to the factor $R_{\bf 1}$, whereas its image should correspond 
to the other factors of the numerator of the $\z$-function. How exactly the geometry of the $\cA$ and $\cB$ curves are related to $\tC_\vph$ will be described elsewhere and will not be pursued in this paper.

We remark further that the map $a$ has as Hodge-component a map 
$$
\a: H^1(\O^2_{\cM_\psi}) \longrightarrow H^0\big(\O^1_{\tC_{\pm\vph}}\big)~.
$$ 
Now the first space 
can be interpreted as the $101$ dimensional space of infinitesimal deformations of quintic $\cM_\psi$, thought
of as the space of degree 5 polynomials $P$ modulo the Jacobian ideal. It follows from the work of H. Clemens that zeros of the holomorphic 1-form $\alpha(P)$ on $\tC_{\pm\vph}$ correspond precisely to the lines that can be infinitesimally lifted over the deformation of $\cM_\psi$ determined by $P$. As the 
curves~$\tC_{\pm\vph}$ both have genus $626$, a differential form has $2{\times}626-2=1250$ zeros. Thus we see that $2{\times}1250=2500$ lines will emerge from the $\tC_\vph$, which together with the $375$ isolated lines gives a total of $2875$ lines that we find on a generic~quintic. 
\newpage
\section{The families of lines}\label{families}
\subsection{Explicit parametrization}\label{Explpar}
Suppose now that, for a line, no coordinate is identically zero. Each $x_i$ is a linear combination of coordinates $(u,v)$ on the line. At least two of the coordinates must be linearly independent as functions of $u$ and $v$. Let us take these independent coordinates to be $x_1$ and $x_2$ then we may take the line to be of the form
\beq
x~=~(u,\, v,\, bu + rv,\, cu+sv,\, du+tv)~.
\label{genline}\eeq
The condition that such a line lies in the quintic imposes the following conditions on the six coefficients:
\beq\begin{split}
 b^5+c^5+d^5+1 ~&=~0\\[3pt]
 b^4 r+c^4 s+d^4 t - b c d\, \psi ~&=~0\\[3pt]
 2\,(b^3 r^2 + c^3 s^2 + d^3 t^2) - (c d r+b d s+b c t)\,\psi ~&=~0\\[3pt]
 2\,(b^2 r^3 + c^2 s^3 + d^2 t^3) - (d r s+b s t+c r t)\,\psi  ~&=~0\\[3pt]
 b r^4+c s^4+d t^4 - r s t\, \psi  ~&=~0\\[3pt]
 r^5+s^5+t^5+1 ~&=~0~.\\ 
\end{split}\label{sixeqs}\eeq
Although there are six equations, we will see that there is a one dimensional family of solutions for the coefficients. However, before coming to this, consider the special case that the coordinates $x_j$ are not all linearly independent as functions of $u$ and $v$. Such a case is equivalent to taking take $r=0$, say, in \eqref{sixeqs}. With this simplification it is straightforward to solve the equations and we find that this case corresponds precisely to the van Geemen~lines.

If we now seek lines that are neither the isolated lines nor the van Geemen lines then we can take all the parameters $b,c,d,r,s,t$ to be nonzero and we also know that all the coordinates are linearly independent as functions of $u$ and $v$. It follows that for a general line, one that is not a isolated line or a van Geemen line, that \eqref{genline} is, in fact, a general form. The first two coordinates of a general line are linearly independent so we choose coordinates so that $x_1=u$ and $x_2=v$ and then the remaining coordinates are linear forms as indicated. Note that we do not have to take separate account of permutations.

In order to simplify \eqref{sixeqs} it is useful to start by scaling the coefficients and the parameter
$$
b~=~c b'~,\quad d~=cd'~,\quad r~=~sr'~,\quad t~=~st'~,\quad \ps~=~cs\ps'~.
$$
This removes $c$ and $s$ from the four central relations. Further scalings lead to additional simplification. This process leads to the following transformation of the variables and parameter
$$
r~=~s\k~,\quad b~=~c\k\t~,\quad d~=~c\k\t\d~,\quad t~=~s\k\t\d\s~,\quad 
\ps~=~\frac{cs}{\d\k^2\t}\,\tilde{\ps}~. 
$$
This has the advantage that, after cancellation, the equations become
\beq\begin{split}
1 + c^5\big[ 1 + \k^5\t^5 (1+\d^5)\big]~&=~0\\[12pt]
1 + \k^5\t^4\, (1\,+\,\d^5\,\s\t\,)  ~&=~\tilde{\ps}\,\t\\[7pt]
1 + \k^5\t^3 (1+\d^5\s^2\t^2)~&=~\frac12\,\tilde{\ps}\,(1+\t+\s\t)\\[7pt]
1 + \k^5\t^2 (1+\d^5\s^3\t^3)~&=~\frac12\,\tilde{\ps}\,(1+\s+\s\t)\\[7pt]
1 + \k^5\t\, (1 + \d^5\,\s^4\t^4)~&=~\tilde{\ps}\,\s\\[10pt]
1 + s^5\big[1 + \k^5(1+\d^5\s^5\t^5)\big]~&=~0~.\\
\end{split}\label{sixeqstransf}\eeq
and depend on $\d$ and $\k$ only through $\d^5$ and $\k^5$. Combining the second, third, fourth and fifth relations with multiples $(1,-2,2,-1)$ results in the cancellation of both the constant and $\tilde{\ps}$ dependent terms. In this way we find
\beq
\d^5~=~\frac{(1-\t)(1-\t+\t^2)}{\s\t^4 (1-\s)(1-\s+\s^2)}~.
\label{delta5}\eeq
Solving the central four relations also for $\k^5$ and $\tilde{\ps}$, we find
\beq
\k^5~=-\frac{(1-\s)(1-\s+\s^2)}{\t (1 - \s\t)(1 - \s\t + \s^2\t^2)}~~~\text{and}~~~
\tilde{\ps}~=~2\,\frac{(1-\s)(1-\t)}{1-\s\t+\s^2\t^2}~.
\label{psitilde}\eeq
Moreover the three relations in \eqref{delta5} and \eqref{psitilde} exhaust the content of the four central equations in \eqref{sixeqstransf}.

The first and last relations in \eqref{sixeqstransf} now give $c$ and $s$ in terms of $\s$ and $\t$. Finally, on substituting what we know into the relation
$$
\ps^5~=~\frac{c^5s^5}{\d^5\k^{10}\t^5}\,\tilde{\ps}^5 ~,
$$
\vskip5pt 
we obtain a constraint $F(\s,\t)=0$ where
\beq\begin{split}
&F(\s,\t) ~=~32\,\s^2\t^2\, (1{-}\s)^2 (1{-}\t)^2 (1{-}\s\t)^2 \,-\\[5pt]
&(1{-}\s{+}\s^2)(1{-}\t{+}\t^2)(1{-}\s\t{+}\s^2\t^2)\!
\Big[1{-}\t(1{+}\s){+}\t^2(1{-}\s{+}\s^2)\Big]\!\!
\Big[1{-}\s(1{+}\t){+}\s^2(1{-}\t{+}\t^2)\Big]\ps^5. \\[7pt]
\end{split}\label{F}\eeq
We are now able to give the lines in terms of $\s$ and $\t$. Let $\a(\s,\t)$ and $\b(\s)$ be given by the~relations
\beq\begin{split}
\a(\s,\t)^5~&=~\s^4\, (1-\s)(1-\t)(1-\s\t)\Big[1-\t(1+\s)+\t^2(1-\s+\s^2)\Big]  \\[5pt]
\b(\s)^5~&=~(1-\s)(1-\s+\s^2)~.
\end{split}\notag\eeq
Then we have
\beq\begin{split}
x_1~&=~\a(\s,\t)\, u \\[5pt]
x_2~&=~\a(\t,\s)\, v \\[5pt]
x_3~&= -\t^\frac45\,\b(\s)\, \left(\s\, u + v\right) \\[5pt]
x_4~&=~~~\b(\s\t)\, \left(\s\, u + \t\, v\right) \\[3pt]
x_5~&= -\s^\frac45\,\b(\t)\, \left(u +\t\, v\right)~. \\
\end{split}\label{fam}\eeq

\Mustata has show that the family of lines has two irreducible components that are isomorphic. This requires $F$ to factorise and this is indeed the case. Setting
\beq
\vph^2~=~\frac{32}{\ps^5} - \frac34
\label{phirelation}\eeq
and
\beq\begin{split}
G~&=~3\s^2\t^2 - \frac12\s\t(1{+}\s)(1{+}\t)(1{+}\s\t) +(1{-}\s{+}\s^2)(1{-}\t{+}\t^2)(1{-}\s\t{+}\s^2\t^2)  \\[5pt]
H~&=~\s\t(1{-}\s)(1{-}\t)(1{-}\s\t)\\
\end{split}\label{GandH}\eeq
we have
$$
F~=-\ps^5\, F_{+}F_{-}~~~\text{with}~~~F_{\pm}~=~G\pm \vph\, H~.
$$

The curves defined by the vanishing of $F_{\pm}(\s,\t)$ are smooth, apart from singularities at the point $(\s,\t)=(1,1)$. Near $(1,1)$ we have the asymptotic form 
\beq 
F_{\pm}(1+\e_1,1+\e_2,\ps)\,\sim\, \e_1^2+\e_1\e_2+\e_2^2~=~(\e_1 - \o\e_2)(\e_1 - \o^2\e_2)~,
\label{odp}\eeq 
so these singularities are ordinary double points. The finite singularities of $F$ are therefore $(1,1)$ together with the solutions of $G=H=0$.

The statement that \eqref{fam} describes all general lines has the following consequence. Clearly if a line can be expressed in the form \eqref{fam} then any permutation of the coordinates $x_k$ yields another line so if the parametrization is general then there must be a reparametrization of $(\s,\t)$ and 
$(u,v)$ that yields this same effect. This is indeed so and the following table gives four such transformations that suffice to generate the permutation group on the~$x_k$. 

The table gives the action of the $\cS_5$ generators on $G$ and $H$. We see that the odd elements of the group interchange $F_{+}$ with $F_{-}$. So each of $F_{\pm}$ is preserved by the alternating 
subgroup~$\cA_5$. Since the odd group elements exchange $F_{+}$ with $F_{-}$ we see that the lines are parametrised by isomorphic curves.

Among the permutations of the $x_k$ there is a cyclic permutation of three coordinates which is of importance. The composition of the exchanges $x_3\leftrightarrow x_5$ and $x_4\leftrightarrow x_5$ generates a cyclic permutation of $(x_3,x_4,x_5)$. As an action on $(\s,\t)$ we have
$$
g_3(\s,\t)~=~\left(\t,\,\frac{1}{\s\t}\right)~.
$$
The action of $g_3$ is expressed most symmetrically by setting $\r=1/\s\t$, so that $\r\s\t=1$, then $g_3$ permutes
$(\s,\t,\r)$ cyclically. We may rewrite the polynomials $G$ and $H$ so as to make the symmetry under $g_3$ manifest. We have
\beq\begin{split}
\frac{G}{(\s\t)^2}~&=~3 - \frac12\,(1+\s)(1+\t)(1+\r) + (1-\s+\s^2)(1-\t+\t^2)(1-\r+\r^2)\\[5pt]
\frac{H}{(\s\t)^2}~&=~ - (1-\s)(1-\t)(1-\r)~.\\
\end{split}\label{P1cubed}\eeq
\subsection{The curves in $\IP^1{\times}\IP^1$ defined by $F_{\pm}$\label{cp1p1}}
We have found curves  in $\IC^2$ defined by $F_{\pm}=0$ whose coverings parametrize lines on  
$\cM_\psi$, with parameters related by \eqref{phirelation}. Let us denote the locus $F_{+}{=}0$ by $C^0_\vph$, the locus $F_{-}{=}0$ is then $C^0_{-\vph}$.

Compactifying $\IC^2$ to $\IP^1{\times}\IP^1$, we obtain a (singular) projective curve of bidegree $(4,4)$. To be explicit, this singular curve is the subset
$$
\left\{\Big((\s_1:\s_2),(\t_1:\t_2)\Big)\,\in\,\IP^1\times\IP^1:\quad 
\s_2^4\t_2^4\, F_{\pm}\left(\frac{\s_1}{\s_2}, \frac{\t_1}{\t_2}\right)~=~0\,\right\}.
$$
The points at infinity are on the lines $\{\infty\}{\times} \IP^1$ and  $\IP^1{\times}\{\infty\}$ (we write $\infty$ for $(1:0)\in\IP^1$):
$$
\begin{array}{ccc}
(\infty,-\o)~,\quad & (\infty,-\o^2)~,\quad & (\infty,0)~,\\[3pt]
(-\o,\infty)~,\quad & (-\o^2,\infty)~,\quad & (0,\infty)~.\\
\end{array}
$$

By means of a Gr\"obner basis calculation, one finds that, for the case that $\cM_\ps$ is smooth, that is 
$\ps^5\neq 1,\infty$, the curves each have three singular points, 
$(\s,\t) = (1,1)$, $(0,\infty)$, $(\infty, 0)$. 
The genus of a smooth bidegree $(d,d')$ curve is $(d-1)(d'-1)$, so if the curve were smooth it would have genus $3{\times} 3=9$. Owing to the singular points, its desingularization has genus at most $6$. The singular points are all related by the operations of \tref{S5transfs} and \eqref{odp} shows the singular points to be ordinary double points, hence the genus of the desingularization is~$9-3=6$. 
Consider now the following list of the 17 points in which the curves $C^0_{\pm\vph}$ intersect (we abuse notation by not distinguishing between the curve in $\IC^2$ and its compactification in $\IP^1{\times}\IP^1$).
$$
\begin{array}{>{}r<{~} >{}r<{~} >{}r<{~} >{}r<{~} >{}r<{~} >{}r<{~} >{}r<{~}}
&&(0,\infty)~, & (\infty,0)~, & (1,\,1)~, \\[3pt]
(0,-\o)~,  & (0,-\o^2)~, & (1,-\o)~, & (1,-\o^2)~, & (-\o,-\o^2)~, & (-\o,\infty)~, & (-\o^2,\infty)~,\\[3pt]
(-\o, 0)~, & (-\o^2,0)~, & (-\o,1)~, & (-\o^2,1)~, & (-\o^2,-\o)~, & (\infty,-\o)~, & (\infty,-\o^2)~.\\
\end{array}
$$ 
\begin{table}
\def\str{\vrule height25pt depth20pt width0pt}
\def\smallstr{\vrule height20pt depth14pt width0pt}
\begin{center}
\begin{tabular}[H]{|>{$\displaystyle}c<{$} |>{$\displaystyle}c<{$} |>{$\displaystyle}c<{$}
|>{$\displaystyle}c<{$}|}
\hline
\multispan4{\smallstr\vrule\hfil\large $\cS_5$ \, generators\hfil\vrule}\\
\hline
\smallstr (\s,\,\t)~~\hbox{transf.} & (u,\,v)~~\hbox{transf.}& \parbox{1.7cm}{\centering effect on coords.}& 
\hbox{effect on $(G,H)$}\\ 
\hline\hline
\str (\t,\, \s) & (v,\, u) &\parbox[c]{1.5cm}{$x_1\leftrightarrow x_2$\\$x_3\leftrightarrow x_5$}
& (G, H)\\ \hline
\str \left(\frac{1}{\s},\,\frac{1}{\t}\right) & (-1)^\frac{1}{5}(\s\t)^\frac{8}{5}\,(v,\, u) 
& x_1\leftrightarrow x_2 
&\frac{1}{\s^4\t^4}\,(G, -H)\\ \hline
\str \left(\frac{1}{\s},\, \s\t\right) & (-\s^\frac{9}{5}\, u ,\, -\s^{-\frac{1}{5}}\, v) & x_4\leftrightarrow x_5 
&\phantom{\t^2} \frac{1}{\s^2}\, (G, -H)\\ \hline
\str \left(\frac{1-\s\t}{1-\t},\, 1-\t\right) 
& \left(\frac{(1-\t)\,(\s u+v)}{(\s\t)^\frac{1}{5}(1-\s\t)^\frac{4}{5}},\; 
-\frac{(1-\s\t)^\frac{1}{5}\,v}{(\s\t)^\frac{1}{5}}\right) & x_1\leftrightarrow x_3
& \left(\frac{\t}{1-\t}\right)^2\!(G,-H)\\ \hline
\end{tabular}
\capt{6.2in}{S5transfs}{The action of four operations, on the coordinates and on the $F_\pm$, that generate 
$\cS_5$.}
\end{center}
\end{table}
\vfill
\begin{table}
\def\str{\vrule height25pt depth20pt width0pt}
\def\smallstr{\vrule height20pt depth14pt width0pt}
\def\nskp{\hskip-5pt}
\begin{center}
\resizebox{6.5in}{!}{
\begin{tabular}[H]{|>{$\displaystyle\nskp}c<{\nskp$} |>{$\displaystyle\nskp}c<{\nskp$} 
|>{$\displaystyle\nskp}c<{\nskp$} |>{$\displaystyle\nskp}c<{\nskp$}|}
\hline
\multispan4{\smallstr\vrule\hfil\large Van Geemen Lines\hfil\vrule}\\
\hline
\smallstr (\s_*,\,\t_*) & \t~~\hbox{for}~~\s=\s_* + \ve & (u, v) & \hbox{Line}\\
\hline\hline
\str (0,\,-\o) & {-}\o {+} 9\g^5\ve
& \left(\frac{c \tilde{u}}{\ve}, -\tilde{v}\right) 
& \big(\tilde{u}, \tilde{v}, -\o^2(c\tilde{u}{-}\tilde{v}), c\tilde{u}{+}\o \tilde{v}, b\tilde{u}\big)\\
\hline
\str (1,\,-\o) &  - \o {+}9\o\g^5\ve 
& \left(-\frac{\o^2\tilde{v}}{\ve^\frac{1}{5}}, - \frac{c\tilde{u}{+}\o\tilde{v}}{\ve^\frac{1}{5}}\right)
&\big(\tilde{v}, c\tilde{u}{+}\o\tilde{v},-\o^2(c\tilde{u}{-}\tilde{v}), b\tilde{u}, \tilde{u}\big)\\
\hline
\str (-\o,-\o^2) &  {-}\o^2 {+} \o \left(\frac{2}{3\ps}\right)^5\!\frac{\ve}{\g^{10}} 
& \left(\frac{c\tilde{u}{+}\o\tilde{v}}{(1{-}\o^2)^\frac{1}{5}\ve^\frac{1}{5}}, 
- \frac{\o^2(c\tilde{u}{-}\tilde{v})}{(1{-}\o^2)^\frac{1}{5}\ve^\frac{1}{5}}\right)
&\big(c\tilde{u}{+}\o\tilde{v},-\o^2(c\tilde{u}{-}\tilde{v}), b\tilde{u},\tilde{v}, \tilde{u}\big)\\
\hline
&\str 1{+}\o^2\ve{+}(\o{+}9\g^5) \ve^2
&\left(\frac{\o c\tilde{u}}{\ve^\frac{6}{5}}{-}\frac{(\o c\tilde{u}{-}\tilde{v})}{2\ve^\frac{1}{5}},
-\frac{\o c\tilde{u}}{\ve^\frac{6}{5}}{-}\frac{(\o c\tilde{u}{-}\tilde{v})}{2\ve^\frac{1}{5}}\right)
&\big(b\tilde{u}, \tilde{u}, \tilde{v}, -\o^2(c\tilde{u}{-}\tilde{v}), c\tilde{u}{+}\o\tilde{v}\big)\\
\smash{\raise27pt\hbox{$(1,\, 1)$}}
&\str 1{+}\o\ve{-}(\o{+}9\g^5)\ve^2
&\left(-\frac{\o c\tilde{u}}{\ve^\frac{6}{5}}{+}\frac{(\o c\tilde{u}{+}\tilde{v})}{2\ve^\frac{1}{5}},
\frac{\o c\tilde{u}}{\ve^\frac{6}{5}}{+}\frac{(\o c\tilde{u}{+}\tilde{v})}{2\ve^\frac{1}{5}}\right)
&\big(\tilde{u}, b\tilde{u}, \tilde{v}, c\tilde{u}{+}\o\tilde{v}, -\o^2(c\tilde{u}{-}\tilde{v})\big)\\
\hline
\end{tabular}
}
\capt{5in}{vanGlines}{The limiting process that gives rise to the van Geemen lines.}
\end{center}
\end{table}
\newpage
\begin{figure}[H]
\begin{center}
\includegraphics[width=6.3in]{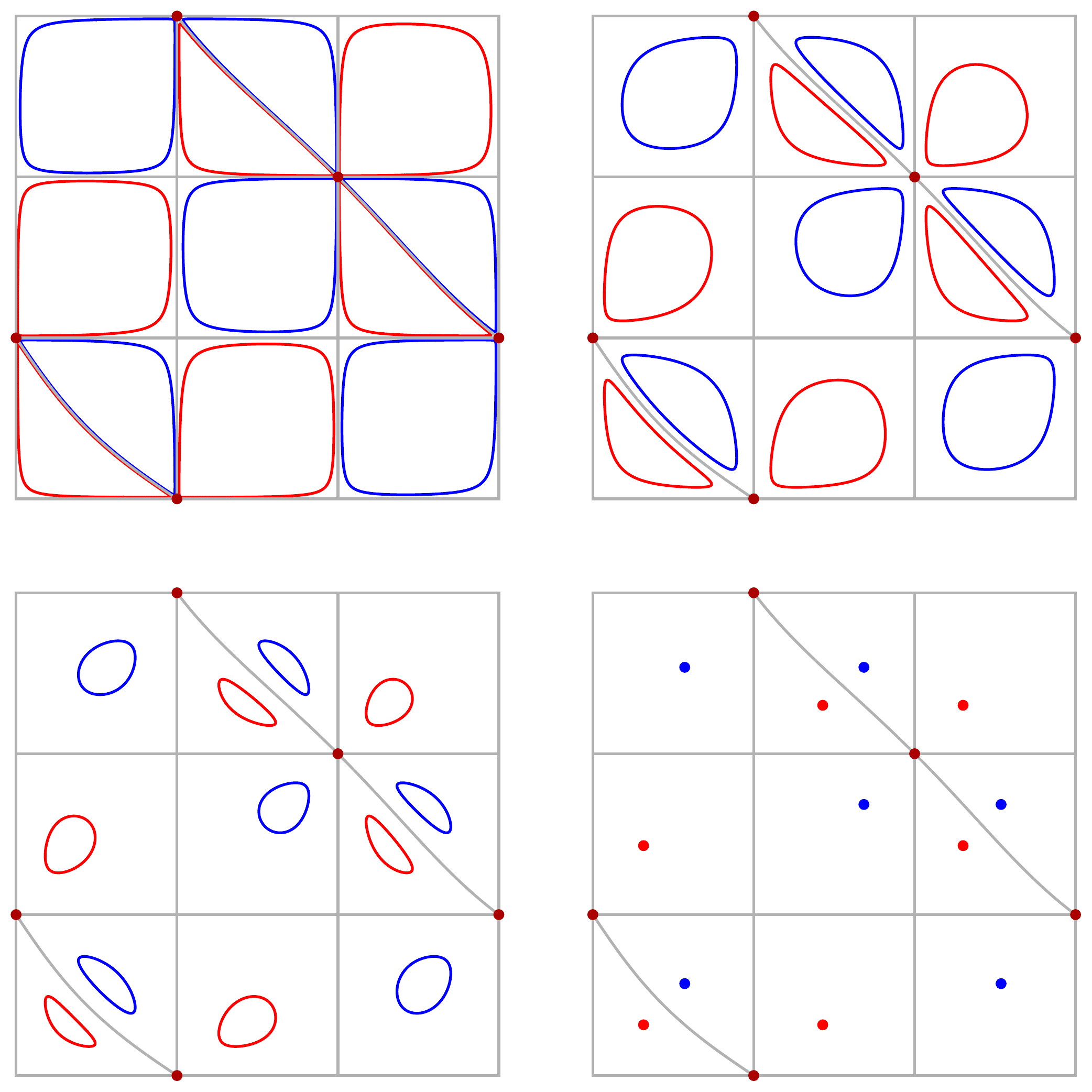}
\vskip0pt
\place{0.2}{3.45}{$\infty$}
\place{1.09}{3.47}{0}
\place{2.03}{3.47}{1}
\place{-0.05}{3.63}{$\infty$}
\place{0.05}{4.55}{0}
\place{0.05}{5.5}{1}
\place{3.55}{3.45}{$\infty$}
\place{4.41}{3.47}{0}
\place{5.35}{3.47}{1}
\place{3.28}{3.63}{$\infty$}
\place{3.37}{4.55}{0}
\place{3.37}{5.5}{1}
\place{0.2}{0.13}{$\infty$}
\place{1.09}{0.13}{0}
\place{2.03}{0.13}{1}
\place{-0.05}{0.3}{$\infty$}
\place{0.05}{1.22}{0}
\place{0.05}{2.17}{1}
\place{3.55}{0.13}{$\infty$}
\place{4.41}{0.13}{0}
\place{5.35}{0.13}{1}
\place{3.28}{0.3}{$\infty$}
\place{3.37}{1.22}{0}
\place{3.37}{2.17}{1}
\capt{6.25in}{animation}{These are plots of the curves $F_{+}=0$, in red, and $F_{-}=0$, in blue, for real 
$(\s,\t)$ as $\ps^5$ ranges from 0 to 1. The diagram is misleading with respect to the points $(1,1)$, $(0,\infty)$ and $(\infty,0)$ which lie on the curve for all $\ps$ but for $\ps\neq 0$ the neighborhoods of the curve on which they lie intersect the plane on which $(\s,\t)$ are both real only in points. The figures show also the images of the 10 exceptional curves of $\dP_5$. These are the 3 points $(1,1)$, $(0,\infty)$ and 
$(\infty,0)$ together with the 7 lines $\s=0,1,\infty$, $\t=0,1,\infty$ and $\s\t=1$. After resolution, the exceptional curves corresponding to the points $(1,1)$, $(0,\infty)$ and $(\infty,0)$ intersect each of the curves $F_\pm =0$ in two points. So too do the other exceptional curves, though the intersections are in complex points not visible in the~figure. The resolved curves are smooth apart from the cases $\ps^5=0,1,\infty$. As 
$\psi\to 0$ the curves tend to the exceptional lines of $\dP_5$ and, as $\psi\to 1$, the curves 
$F_\pm =0$ each develop 6 nodes corresponding to the limiting points shown in the final~figure.}
\end{center}
\end{figure}
The list of the points of intersection includes the three points, just discussed, in which the curves are both singular. Note that these points do not depend on $\vph$.

We know that at least some of the van Geemen lines must lie in the continuous families. Indeed, \Mustata has shown that they all lie in the continuous families, since the only isolated lines are the 375 lines that we have identified as such.
The van Geemen lines are, however, not easy to see from the parametrisation \eqref{fam}. It is a surprising fact that these lines appear precisely as limits, as we approach the points in which the curves $C^0_{\pm\vph}$ intersect. For the points $(0,-\o)$, $(1,-\o)$, $(-\o,-\o^2)$ and the singular point $(1,1)$, this resolution is given in \tref{vanGlines}. All the other resolutions may be obtained from these by acting with the $\cS_5$ operations of \tref{S5transfs}. Each of the nonsingular points of intersection $(\s_*,\t_*)$ gives rise to two van Geemen lines, one in each of the families. The two possible values 
$$
\g^5~=~\frac19\left( \frac12\mp\frac{\ii\vph}{\sqrt{3}} \right)
$$
correspond, respectively, to the two curves $C^0_{\pm\vph}$. For the three singular points each curve has self intersection so the resolution produces two lines for each curve, again the two choices for $\g^5$, as above, correspond, respectively, to the two curves $C^0_{\pm\vph}$. In this way we find 
$14{\times}2{+}3{\times}4=40$ lines which become $40{\times}125=5000$ lines under the action of $\cG$. Thus we have found all the van Geemen lines as resolutions of intersection of the curves $C^0_{\pm\vph}$.

The appearance of fifth roots in \eqref{fam} indicates that we have to allow for different branches and the effect of fifth roots of unity. In \eqref{fam} we have to choose a fifth root of unity for each of $\s$, $\t$, $\a(\s,\t)$ and $\b(\s)$. This might suggest a $\IZ_5^4$ covering, however multiplying all the coordinates $x_j$ by a common factor is of no consequence, so there is in fact a $\IZ_5^3$ covering and we
can allow for different branches of solutions by acting with $\cG$ on a given branch. Somewhat surprisingly monodromy around the singularities of $F_{\pm}$ does not generate~$\cG$. Instead the monodromy simply multiplies all the components $x_j$ of a line by a common factor of $\z^k$ for some $k$. Thus there is no local ramification of the solution. We will give a better description of the 125:1 cover in \SS\ref{cover125}.
\subsection{A partial resolution of the singularities of $C^0_\vph$}
We have seen that the curves $C^0_\vph$ have three singular points. We wish to resolve these singularities. It is interesting to note that two of these singularities can be resolved very naturally. It was remarked previously that by introducing a new parameter $\r$, subject to the constraint $\r\s\t=1$, the equations $F_\pm=0$ can be written, as in \eqref{P1cubed}, so as to be manifestly symmetric under an 
$\cS_3$ subgroup of the permutation symmetry. Once we introduce $\r$, we are dealing with the nonsingular surface $\r\s\t=1$ embedded in $(\IP^1)^3$. If written in homogeneous coordinates, this surface is given by the trilinear equation
\beq
\s_1\t_1\r_1~=~\s_2\t_2\r_2~.
\label{TrilinearEq}\eeq
The vanishing locus of a nonsingular trilinear polynomial in $(\IP^1)^3$ is isomorphic to the del~Pezzo surface 
$\dP_6$, which we may think of as $\IP^2$ blown up in three points. Two of these blow ups resolve the singularities at $(\s,\t)=(0,\infty)$ and $(\s,\t)=(\infty,0)$. Consider the first of these singularities. In homogeneous coordinates the location of the singularity is
$$
\Big((\s_1,\s_2),\,(\t_1,\t_2)\Big)~=~\Big((0,1),\,(1,0)\Big)~.
$$
For these values \eqref{TrilinearEq} is satisfied for all values of $(\r_1,\r_2)$, so the singular point has been replaced by an entire $\IP^1$. A Gr\"obner basis calculation shows that the curves defined by $F_\pm=0$ are now only singular at the point $(\s,\t,\r)=(1,1,1)$.

The surfaces $\dP_6$, $\IP^1{\times}\IP^1$ and $\IP^2$ are all toric and it is clear from their respective fans that $\dP_6$ is obtained from $\IP^2$ by blowing up three points and may also be obtained from 
$\IP^1{\times}\IP^1$ by blowing up two points (for the relation between the blow ups of $\IP^2$ and 
$\IP^1{\times}\IP^1$ see \SS\ref{P2blowup}). Since we wish to resolve the remaining singularity of the curves $F_\pm =0$, it is natural to blow up one further point. This brings us to a consideration of $\dP_5$.
\newpage
\section{The quintic del Pezzo surface $\dP_5$}\label{dp5}
\subsection{Blowing up three points in $\IP^1{\times}\IP^1$}
The curves $C^0_\vph$ in $\IC^2$ define singular curves of bidegree $(4,4)$ in $\IP^1{\times}\IP^1$ which in general have three ordinary double points in $(\sigma,\tau)=(1,1),(0,\infty),(\infty,0)$. The blow up of $\IP^1{\times}\IP^1$ in these three points is the quintic del Pezzo surface $\dP_5$. 

The blow up is given by the polynomials of bidegree $(2,2)$ which are zero in these three points (see Section \ref{picdp5}).
The polynomials of bidegree $(2,2)$ are a $9=3^2$-dimensional vector space with basis $\sigma_1^{a}\sigma_2^b\tau_1^c\tau_2^d$, $a+b=2=c+d$.
The blow up map can thus be given by 
$$ 
\Psi:\,\IP^1\times\IP^1\,\dashrightarrow\, \dP_5\quad(\subset\IP^5)~,\qquad 
(\sigma,\tau)\,\longmapsto\,(z_0,\ldots,z_5)~,
$$
with the $6$ functions (written inhomogeneously for simplicity):
$$ z_0 := \sigma^2\tau^2 - 1,\quad
   z_1 :=  \sigma\tau^2 - 1,\quad
   z_2 :=  \sigma^2\tau - 1,\quad
    z_3 := \sigma\tau - 1,\quad
   z_4 :=  \tau - 1,\quad
   z_5 := \sigma - 1.
$$
The image of $\IP^1{\times}\IP^1$ is $\dP_5$, in its anti-canonical embedding into $\IP^5$.
To find the inverse, notice that $(z_3-z_5,z_4)=(\sigma\tau-\sigma,\tau-1)=(\sigma,1)$ in $\IP^1$.
Thus the inverse map $\Phi$, which is everywhere defined, is given by
$$
\Phi:\,\dP_5\,\longrightarrow\,\IP^1\times\IP^1~,\qquad z:=(z_0,\ldots,z_5)\,\longmapsto\,
\Big((z_3-z_5,z_4),\,(z_3-z_4,z_5)\Big)
$$
(one should notice however that this formula for $\Phi$ only works on an open subset of $\dP_5$, using certain quadratic relations between the $z_i$ which are satisfied on $\dP_5$, one can extend $\Phi$ to all of $\dP_5$).

The surface $\dP_5\subset\IP^5$ is defined by 5 quadratic equations.
An example of such an equation is
$$
q_0\,=\,0~,\qquad\mbox{with}\quad q_0\,:=\,(z_1-z_3)z_5\,-\,(z_2-z_3)z_4~,
$$
in fact $\Big((\s\t^2-1)-(\s\t-1)\Big)(\s-1)=\Big((\s^2\t-1)-(\s\t-1)\Big)(\t-1)$.

The image of the curve $C^0_\vph$ is defined by an additional quadratic equation,
which we will discuss in section \ref{eqndp}.

\subsection{Automorphisms of $\dP_5$}\label{autdp}
As we will recall below, the group of automorphisms of the algebraic surface $\dP_5$ is $\cS_5$.
The action of $\cS_5$ on $\IP^1\times\IP^1$, given by the birational transformations given in the previous section,  induces these automorphisms on $\dP_5$.
The automorphisms of $\dP_5$ act linearly on the $z_i$'s 
(since they are the sections of the anti-canonical bundle of $\dP_5$). Thus we get a much simpler description of the $\cS_5$ action.

We will now determine the matrices of the four elements in $\cS_5$ given in \tref{S5transfs}.
One should notice that, for example, $(z_0,\ldots,z_5)$ and $(-z_0,\ldots,-z_5)$ define the same point in $\IP^5$, but they are distinct as points in $\IC^6$. To obtain a linear representation of $\cS_5$ on $\IC^6$ one has to make the choices we give below.

The element $(12)(35)$ acts as $(\sigma,\tau)\mapsto(\tau,\sigma)$ on $\IP^1{\times}\IP^1$
and as
\begin{align*}
(12)(35):\quad z\,&\longmapsto\,(-z_0,\,-z_2,\,-z_1,\,-z_3,\,-z_5,\,-z_4)\\
\intertext{on $\IC^6$. Notice that the trace of $(12)(35)$ on $\IC^6$ is $-2$.
The second permutation is $(12)$ which acts as $(\sigma,\tau)\mapsto(\sigma^{-1},\tau^{-1})$, so as 
$\big((\sigma_1,\sigma_2),\,(\tau_1,\tau_2)\big)\mapsto \big((\sigma_2,\sigma_1),\,(\tau_2,\tau_1)\big)$ in homogeneous coordinates. This gives the map, with trace zero,}
(12):\quad z\,&\longmapsto\, (-z_0,\,-z_0+z_5,\,-z_0+z_4,\,-z_0+z_3,\,-z_0+z_2,\,-z_0+z_1)~.\\
\intertext{The permutation $(45)$ acts non-linearly, $(\sigma,\tau)\mapsto(\sigma^{-1},\sigma\tau)$,
substituting this in the polynomials $z_i$ and multiplying by $-\sigma$ gives the action on $\IC^6$:}
(45):\quad z\,&\longmapsto\, (-z_1+z_5,\,-z_0+z_5,\,-z_4+z_5,\,-z_3+z_5,\,-z_2+z_5,\,z_5)~.\\
\intertext{Finally we have $(13)$ acting as $(\sigma,\tau)\mapsto \big((1-\sigma\tau)(1-\tau),1-\tau\big)$, substituting and multiplying by $(1-\tau)/\tau$ gives the linear map:}
(13):\quad z\,&\longmapsto \,
(-z_0+z_2+2z_3-2z_5,\, -z_1+2z_3+z_4-z_5,\,z_2-2z_5,\,z_3-z_5,\,z_4,-z_5)~.\\[-15pt]
\end{align*}
We have verified that this gives indeed a linear representation of $\cS_5$ on $\IC^6$. Computing the traces and comparing with a character table of $\cS_5$ (see section \ref{eqndp}), we find that this representation is the unique irreducible $6$-dimensional representation of $\cS_6$. 

\subsection{Exceptional curves in $\dP_5$}
We obtained $\dP_5$ as the blow up of $\IP^1\times\IP^1$ in the three points $(1,1),(0,\infty),(\infty,0)$. Thus on $\dP_5$ we have three $\IP^1$'s, the exceptional curves over these points. These are lines in $\IP^5$ lying on $\dP_5$. To find them, it suffices to find just one and then apply suitable elements of $\cS_5$ to find the others. 
The points $(a:b)$ on the exceptional curve over $(1,1)$ are the limit points of the image of $(\sigma,\tau)=(1+\epsilon a,1+\epsilon b)$ for $\epsilon\rightarrow 0$ under the blow up map. One finds the line
$$
E_{12}\,:\quad (2a+2b,\,a+2b,\,2a+b,\,a+b,\,b,\,a),\qquad (a,\,b)\,\in\,\IP^1~.
$$
In fact, $\Psi(E_{12})=\big((a+b-a,b),\,(a+b-b,a)\big)=\big((b,b),\,(a,a)\big)=\big((1,1),(1,1)\big)$ which is indeed $(1,1)$ in inhomogeneous coordinates.
From equation \eqref{odp} we infer that the (strict transforms of the) curves $C^0_\vph$ intersect $E_{12}$ in two points, independent of $\vph$, which correspond to 
$(a:b)=\big((\omega,1),(\omega^2,1)\big)$. In the following we shall give parametrisations of the other exceptional curves. In each case, the parameters $(a,b)$ will be understood as the coordinates of a $\IP^1$.
\begin{table}[t]
\def\bigstr{\vrule height18pt depth12pt width0pt}
\def\str{\vrule height15pt depth8pt width0pt}
\begin{center}
\begin{tabular}[H]{| >{$}c<{$} | >{$}c<{$} | >{$}c<{$} | >{$}c<{$} |}
\hline
\multispan4{\bigstr\vrule\hfil\large  Exceptional curves in $\dP_5$ \hfil\vrule}\\
\hline
\str \hbox{Name} & \hbox{Parametrization} & \hbox{Image in}~\IP^1\times \IP^1 & \hbox{Special points}\\
\hline\hline
\str E_{12} &(2a+2b,a+2b,2a+b,a+b,b,a)
& (1,1)
& \hbox{singular point}\\
\hline
\str E_{13} &  (a,b,0,0,0,0) 
& \t\,=\,\infty
&(-\o,\infty),\,(-\o^2,\infty) \\
\hline
\str E_{14} &  (0,0,a,0,0,b)
& (\infty,0)
& \hbox{singular point}\\
\hline
\str E_{15} & (a,a,a,a,b,a)
& \s\,=\,0
& (0,-\o),\,(0,-\o^2)\\
\hline
\str E_{23} & (a,a,a,a,a,b)
&  \t\,=\,0
& (-\o,0),\,(-\o^2,0) \\
\hline
\str E_{24} & (0,a,0,0,b,0)
& (0,\infty)
&\hbox{singular point} \\
\hline
\str E_{25} & (a,0,b,0,0,0)
& \s\,=\,\infty
& (\infty,-\o),\,(\infty,-\o^2) \\
\hline
\str E_{34} & (a,\,b,\,a,\,b,\,0,\,b)
& \tau\,-\,1\,=\,0
&(-\o,1),\,(-\o^2,1)\\
\hline
\str E_{35} & (0,\,b,\,a,\,0,\,b,\,a)
& \s\t\,-\,1\,=\,0
&  (-\o,-\o^2),\,(-\o^2,-\o)\\
\hline
\str E_{45} & (a,a,b,b,b,0)
& \s\,-\,1\,=\,0
& (1,-\o),\,(1,-\o^2)\\
\hline
\end{tabular}
\capt{5.7in}{ExcCurves}{The ten exceptional curves in $\dP_5$, showing their images in $\IP^5$ and in $\IP^1{\times}\IP^1$. The table also gives the points in which the divisors meet the curve $C_\vph^0$.}
\end{center}
\end{table}

One verifies that this line is mapped into itself under the action of $(12),(34),(45)\in \cS_5$, which generate a subgroup of order $2{\times}6=12$ in $\cS_5$. Acting with elements of $\cS_5$ on $E_{12}$ produces $9$ other lines, which are denoted by $E_{ij}=E_{ji}$, $1\leq i,j\leq 5$ and $i\neq j$, compatible with the action of $\cS_5$. 

We now discuss some of these lines in $\dP_5$ and their source in $\IP^1\times \IP^1$.
The line in $\dP_5$ which is the exceptional curve over $(0,\infty)$ can be found with a limit as above and it is
\begin{align*}
E_{24}\,&:\quad (0,\,a,\,0,\,0,\,b,\,0)~,\\
\intertext{again one verifies easily that $\Phi(E_{24})=\big((0-0,b)),\,(0-b,0)\big)=\big((0,1),\,(1,0)\big)$ which is~$(0,\infty)$. As $(12)(35)$ permutes $\sigma$ and $\tau$, and thus $(0,\infty)$ and $(\infty,0)$, the exceptional curve over $(\infty,0)$~is}
E_{14}\,&:\quad (0,\,0,\,a,\,0,\,0,\,0,\,b)~.\\
\end{align*}
\vskip-10pt
The rulings $(1,1){\times}\IP^1$ and $\IP^1{\times}(1,1)$ passing through $(1,1)$ are also mapped to lines, for example, in inhomogeneous coordinates:
$$
\Psi(1,\,a)\,=\,(a^2-1,\,a^2-1,\,a-1,\,a-1,\,a-1,\,0)\,=\,(a+1,\,a+1,\,1,\,1,\,1,\,0)~,
$$
which shows that the curve defined by $\sigma=1$ maps to a line, which is $E_{45}$, on $\dP_5$:
\begin{align*}
E_{45}\,&:\quad (a+b,\,a+b,\,b,\,b,\,b,\,0)~.\\
\intertext{Similarly, the curve $\tau=1$ (obtained from the first by $(12)(35)\in\cS_5$) maps to the line}
E_{34}\,&:\quad (a+b,\,b,\,a+b,\,b,\,0,\,b)~.\\
\intertext{In this way each of the three points $(1,1),(0,\infty),(\infty,0)$ provides us with three lines on 
$\dP_5$, so we already have $9$ lines. For example, the curve $\tau=\infty$ maps to the line}
E_{13}\,&:\quad (a,b,0,0,0,0)~.\\
\intertext{A final $10$th line is given by the image of the unique curve of bidegree $(1,1)$ passing through these three points. Its equation is $\s_1\t_1- \s_2\t_2=0$, i.e.\ $\sigma\tau=1$,
so it can be parametrized by $(a,a^{-1})$ and its image under $\Psi$ is
$$
\Psi(a,a^{-1})\,=\,(0,\,a^{-1}-1,\,a-1,\,0,\,a^{-1}-1,\,a-1)\,=\,
(0,\,1,\,-a,\,0,\,1,\,-a)~.
$$
Thus we have found the line}
E_{35}\,&:\quad (0,\,b,\,-a,\,0,\,b,\,-a)~.\\
\end{align*}
\subsection{The curves $C_\vph$ and the Wiman pencil}\label{eqndp}
We will use representation theory of $\cS_5$ to find the equations of the curves $C_\vph$.

The coordinates on $\IP^5$ are $z_0,\ldots,z_5$ and the action of $\cS_5$ on these coordinates was given in section \ref{autdp}. Comparing the traces with \tref{CharTab}, we find that the linear functions are in the unique $6$-dimensional irreducible representation of $\cS_5$.
\begin{table}
\begin{center}
\def\bigstr{\vrule height17pt depth10pt width0pt}
\def\str{\vrule height13pt depth6pt width0pt}
\begin{tabular}[H]{|>{\hskip5pt $\bf}l <{$}|>{\hskip5pt}c<{\hskip5pt}|>{\hskip3pt} r<{\hskip8pt}|r<{\hskip15pt}
| >{\hskip5pt}r<{\hskip13pt}|r<{\hskip12pt}| r<{\hskip15pt}| r<{\hskip20pt}|}
\hline
\multispan8{\bigstr\vrule\hfil\large Characters of $\cS_5$\hfil\vrule}\\ \hline
\str&$e$&(12)\hspace*{-5pt}&(12)(34)\hspace*{-15pt}&(123)\hspace*{-10pt}&(1234)\hspace*{-12pt}&(12345)\hspace*{-15pt}&(123)(45)\hspace*{-20pt}\\
\hline\hline
\str 1     &1&1&1&1&1&1&1\\ \hline
\str 1_b  &1&-1&1&1&-1&1&-1\\ \hline
\str 4     &4&2&0&1&0&-1&-1\\ \hline
\str 4_b &4&-2&0&1&0&-1&1\\ \hline
\str 5     &5 &1 &1&-1& -1&  0&  1\\ \hline
\str 5_b &5&-1 &1&-1 &1 & 0& -1\\ \hline
\str 6     &6&0&-2&0&0&1&0\\ \hline
\end{tabular}
\vskip8pt
\capt{4.5in}{CharTab}{The character table of $\cS_5$. This proves useful in identifying the image of $C_\vph$ in $\IP^5$.}
\end{center}
\end{table}
The 21-dimensional representation $S_2$ of $\cS_5$ on the quadratic functions $z_iz_j$ 
has character $\chi_2$ given by $\chi_2(g)=(\chi(g)^2+\chi(g))/2$, 
where $\chi$ is the character of $\cS_5$ on the linear functions.
Decomposing it into irreducible characters one finds:
$$
S_2\,=\,{\bf 1}\oplus{\bf 1_b}\oplus{\bf 4}\oplus2\cdot{\bf 5}\oplus{\bf 5_b}~. 
$$

Let $G_z,H_z\in I_2$ be polynomials which span $\bf 1$ and $\bf 1_b$ respectively.
Thus $G_z$ is $\cS_5$-invariant and hence $G_z=0$ defines a $\cS_5$ curve in $\cP_5$.
Similarly, as $gH_z=\epsilon(g) H_z$, 
where $\epsilon(g)$ is the sign of the permutation $g$,
the curve $H_z=0$ is $\cS_5$-invariant. 
Such polynomials can be found as $\sum_g g(z_0z_1)$ and
$\sum_g \epsilon(g)g(z_0z_1)$, where the sum is over all $g\in \cS_5$.
To relate these polynomials in the $z_i$ to those in $\s,\t$, recall that
the $z_i$ correspond to a basis of the polynomials of bidegree $(2,2)$ 
which vanish in $(1,1),(0,\infty),(\infty,0)$. More precisely, using the map $\Psi$,
we have
$$
\Psi^*(z_0)\,=\,\s^2\t^2-1,\qquad\ldots,\quad \Psi^*(z_5)\,=\,\s-1~.
$$

The unique $\cS_5$ invariant quadratic polynomial $G_z$ is:
\beq\begin{split}
G_z~:=~2z_0^2 &- 2z_0z_1 - 2z_0z_2 - 2z_0z_3 + z_0z_4 + z_0z_5 + 2z_1^2 + z_1z_2 - 2z_1z_3 \\
& - 2z_1z_4 + 2z_2^2 - 2z_2z_3 - 2z_2z_5 + 6z_3^2 - 2z_3z_4 - 2z_3z_5 + 2z_4^2 + z_4z_5 + 2z_5^2~,\\
\end{split}\notag\eeq
and we have verified that 
$$
\Psi^*G_z\,:=\,G_z(\s^2\t^2-1,\ldots,\s-1)\,=\,G(\s,\t)~.
$$
Similarly, the unique quadratic polynomial $H_z$ invariant under $\cS_5$, up to a sign, is
$$
H_z~:=~\frac13(-2z_0z_3 + z_0z_4 + z_0z_5 - z_1z_2 + 2z_1z_3 + 2z_2z_3 - 2z_3z_4 - 2z_3z_5 + z_4z_5)
~,
$$
and one finds that
$$
\Psi^*G_z\,=\,H(\s,\t)~,
$$
where, in the above, $G(\s,\t)$ and $H(\s,\t)$ are the polynomials \eqref{GandH}.

The curves $C^0_\vph$ in $\IP^1\times\IP^1$ which have equation 
$F_+=G+\vph H=0$
are thus the images under the blow down 
$\Phi:\dP_5\rightarrow \IP^1\times\IP^1$ of the curves 
$C_\vph$ in $\dP_5$ defined by $G_z+\vph H_z=0$. 
This pencil of curves $\{C_\vph\}_{\vph\in\IP^1}$ is known as the Wiman pencil. 
The curve $C_0$, defined by $G_z=0$, is smooth and has automorphism group $\cS_5$, 
and is known as the Wiman curve. 

The curves $C_\vph$ have a 125:1 cover $\tC_\vph$ which parametrizes lines on the
quintic threefold $\cM_\psi$ Dwork pencil, where $\vph$ and $\psi$ are related as in Section \ref{Explpar}.
We will turn to this covering in section \ref{cover125}. 

We conclude with one final remark on the 
curve $C_\infty$ defined by $H_z=0$ in $\dP_5$. The homogeneous polynomial defined by $H$ of must be of bidegree $(4,4)$, 
so besides the 5 factors in the dehomogenized equation, we should take into account two more: 
$$
H(\s_1,\s_2,\t_1,\t_2)\,=\,\s_1\t_1(\t_2-\s_1)(s_2-\t_1)(\s_2\t_2-\s_1\t_1)\s_2\t_2
$$
Thus the curve $H=0$ actually has $7$ irreducible components which all map to lines in $\dP_5$ as we observed earlier. Moreover, $H=0$ passes through the three points which get blown up. 
In fact one can check (see also \sref{picdp5}) that the curve $H_z=0$ in $\dP_5$ has 10 irreducible components, which are the 10 lines in $\dP_5$, each with multiplicity one.

On each of the $10$ lines in $\dP_5$ there are two points which correspond to van Geemen  
lines. Each line is invariant under a subgroup of order $12$ of $\cS_5$ and these two points are the fixed point set of any of the two elements of order three in the subgroup.

\newpage
\section{A second look at the curves parametrizing the lines}\label{secondlook}
\subsection{The Pl\"ucker map}\label{plmap}
The explicit parametrization of the lines in the Dwork pencil was given in Eq~\eqref{fam}.
We will now study their Pl\"ucker coordinates, which will be the key to understanding the 125:1 cover $\tC_\vph\rightarrow C_\vph$.

Given a line $l$ in $\IP^4$ spanned by two points 
$x=(x_1,\ldots,x_5)$ and $y=(y_1,\ldots,y_5)$,
its Pl\"ucker coordinates $\pi_{ij}(l)=-\pi_{ji}(l)$ are defined as:
$$
\pi_{ij}(l)\,:=\,x_iy_j\,-\,x_jy_i,\qquad l=\langle x,y\rangle\quad\subset \IP^4.
$$
The $10$ Pl\"ucker coordinates $\pi_{ij}(l)$ with $1\leq i<j\leq 5$, 
viewed as projective coordinates on $\IP^9$, determine $l$  uniquely. 

The van Geemen lines given in equation (\ref{VanGone}) are spanned by the rows of the matrix (the first corresponding to the point with $(u,v)=(1,0)$ and the second to $(u,v)=(0,1)$)
$$
\left(
\begin{array}{ l<{\hskip15pt} l >{\hskip5pt}c<{\hskip5pt} l<{\hskip5pt} l } 
1 & 0 &\z^{-k-l}b &\z^k c &-\z^l \o^2 c \\[5pt]
0 & 1 & 0              &\z^k\o &\phantom{-}\z^l\o^2
\end{array}
\right)\!\raisebox{-14pt}{.}
$$
One notices that $\pi_{13}=0$, and that the other $\pi_{ij}$ are non-zero. 
Recall that these lines are invariant under the cyclic permutation of 
$(x_2,x_4,x_5)$. 
As the other van~Geemen lines are obtained from this one by the action of the group $\cS_5{\rtimes}\cG$, we conclude that a van~Geemen line has exactly one Pl\"ucker coordinate $\pi_{ij}$ which is zero, the indices $ij$ are such that the stabilizer of the line is conjugated in $\cS_5{\rtimes}\cG$ to the cyclic subgroup generated by $(klm)\in\cS_5$ where $\{i,j,k,l,m\}=\{1,\ldots,5\}$.

These indices $i,j$ can also be obtained as follows.  
The point on $C_\vph$ determined by such a line lies on the intersection of this curve with the line $E_{pq}$ on $\dP_5$ and the sets of indices $\{i,j\}$ and $\{p,q\}$ are the same. We will now see that, conversely, a line in the Dwork pencil for which one of the $\pi_{ij}$ is zero are is a van Geemen line.

The elements of the group $\cG$ acts by multiplying the coordinates 
$x_1,\ldots,x_5$ of $\IP^4$ by fifth roots of unity. Hence the induced action
of $\cG$ on the Pl\"ucker coordinates is also by multiplication by fifths roots of unity.
The fifth powers $\pi_{ij}^5$ of the Pl\"ucker coordinates are thus invariant under $\cG$ and hence functions on $C_\vph$, more precisely, the quotients $\pi_{ij}^5/\pi^5_{pq}$ define meromorphic functions on $C_\vph$. 
These functions are easy to find.

The Pl\"ucker coordinates of the lines parametrized by the 125:1 cover of $C_\vph$ in \eqref{fam} are the determinants of the $2\times 2$ minors of the following matrix: 
$$
\left(
\begin{array}{ c c >{\hskip3pt} l >{\hskip7pt}l >{\hskip3pt}l}
\alpha(\s,\t) & 0                    & -\t^{4/5}\beta(\s)\s &\beta(\s\t)\s  &-\s^{4/5}\beta(\t)\\[12pt]
0                 &\;\alpha(\t,\s)\;&-\t^{4/5}\beta(\s)     &\beta(\s\t)\t  &-\s^{4/5}\beta(\t)\t \\
\end{array}
\right)\!\raisebox{-18pt}{.}
$$

From this we compute for example, with $\beta(\s)^5=(1-\s)(1-\s+\s^2)$:
$$
\pi_{35}^5\,=\,\Big(\t^{4/5}\beta(\s)\s{\cdot}\s^{4/5}\beta(\t)\t - 
\s^{4/5}\beta(\t){\cdot}\t^{4/5}\beta(\s)\Big)^5
\,=\,\s^4\t^4\beta(\s)^5\beta(\t)^5(\s\t-1)^5~.
$$
In this way we get $10$ polynomials in $\s,\t$ of rather high degree, 
but they do have a common factor, which is:
$$
p_c\,:=\,\s^4\t^4(\s-1)(\t-1)(\s\t-1)~.
$$
The quotients $\pi_{ij}^5/p_c$ can be homogenized to polynomials 
of bidegree $(6,6)$ in $\s_1,\s_2$ and $\t_1,\t_2$:
$$
p_{ij}(\s_1,\s_2,\t_1,\t_2)\,:=\,(\s_2\t_2)^6(\pi_{ij}^5/p_c)(\s_1/\s_2,\t_1/\t_2)~.
$$
These $p_{ij}$ are reducible. Their irreducible components can be used to define meromorphic functions on $C_\vph$ 
with quite interesting zeroes and poles as we will see in the next sections and in the Appendix. 

We introduce some notation for the irreducible components of the $p_{ij}$.
The polynomial defining the curve in $\IP^1{\times}\IP^1$ which maps to the line 
$E_{ij}$ on $\dP_5$ is denoted by $m_{ij}$, and we give them in \tref{Divs}.
We have two polynomials, of bidegree $(2,0)$ and $(0,2)$ respectively,
$$
l_1\,:=\,\s_1^2 - \s_1\s_2 + \s_2^2,\qquad
l_2\,:=\,\t_1^2 - \t_1\t_2 +\t_2^2
$$
which are reducible (\,$l_1=(\s_1+\o\s_2)(\s_1+\o^2\s_2)$\,) and $l_i=0$ intersects the curve $C_\vph^0$ in special points corresponding to van Geemen lines. 

The intersection of $C_\vph$ with another curve is written as a divisor $\sum_P n_pP$, which is a formal finite sum with $P\in C_\vph$ and $n_p\in\IZ$ the intersection multiplicity. We write $D_{ij}=P_{ij}+Q_{ij}$
for the divisor given by the pair of special points on $C_\vph$ which correspond to the van Geemen lines indexed by $ij$. Thus if $ij=45$, we can take $P_{ij}=(1,-\o)$ and $Q_{ij}=(1,-\o^2)$, viewed as points on the smooth model $C_\vph$ of $C^0_\vph$.
In case $ij=12$ we take the two points of $C_\vph$ which map to the singular point $(1,1)$ of $C^0_\vph$. On $\dP_5$, these divisors are the intersection divisors of $C_\vph$ and the lines $E_{ij}$:
$$
D_{ij}\,:=\,C_\vph\,\cap\,E_{ij}~.
$$
However, pulling back the divisors $m_{ij}=0$, we do not get the divisors $E_{ij}$, but we also get contributions from the singular points.
In Table \ref{Divs} we give the precise results.

With this notation, table \ref{ExcCurves} shows that
$$
(l_1=0)\,\cap\,C_\vph \,=\,D_{13}\,+\,D_{23}\,+\,D_{34}\,+\,D_{35}~,
$$
applying $(12)(35)\in\cS_5$ one obtains $(l_2=0)\cap\,C_\vph$.

Finally there are three polynomials, of bidegree $(2,2)$, which turn out to be reducible.
The first is
{\renewcommand{\arraystretch}{1.5}
$$
\begin{array}{rcl}
k_{14}&:=&\s_1^2\t_1^2 - \s_1^2\t_1\t_2 + \s_1^2\t_2^2 - \s_1\s_2\t_1\t_2 - \s_1\s_2\t_2^2 + \s_2^2\t_2^2\\
&=&
(\s_1\t_1 + \o^2\s_1\t_2 + \o\s_2\t_2)
(\s_1\t_1 + \o\s_1\t_2 + \o^2\s_2\t_2)~.
\end{array}
$$
}
The first factor defines a curve in $\IP^1{\times}\IP^1$ which can be parametrized by
$$
(s,t)\,\longmapsto\,\Big((\s_1,\s_2),\,(\t_1,\t_2)\Big)
\,=\, \Big((-\o t ,\o^2t+s),\,(s,t)\Big)~,\qquad(s,t)\,\in\,\IP^1~.
$$
The intersection of this curve with $C^0_\vph$, which is defined by $F_+=0$, is obtained from
$$
F_+\big((-\o t ,\o^2t+s),\,(s,t)\big)\,=\,(2\vph - 2\o - 1)st^3(s-t)^3(s +\o^2t)~.
$$
One finds the points $(-\o^2,0)$, $(\infty,-\o^2)$, which are in the divisors $D_{23}$ and $D_{25}$ respectively with multiplicity one and the singular points
$(0,\infty)$ and $(1,1)$ with multiplicity three. So the curve must be tangent to one branch of $C^0_\vph$ in these points. The equation of the other factor is the complex conjugate, so we conclude that
$$
(k_{14}=0)\,\cap\,C_\vph \,=\,D_{23}\,+\,D_{25}\,+\,3D_{24}+3D_{12}~.
$$
\begin{table}
\def\bigstr{\vrule height18pt depth12pt width0pt}
\def\str{\vrule height14pt depth9pt width0pt}
\begin{center}
\begin{tabular}[H]{|>{$}c<{$} | >{$~}l<{~$} |>{$}l<{$} |}
\hline
\multispan3{\bigstr\vrule\hfil\large Curves and intersection divisors with $C_\vph$ \hfil\vrule}\\
\hline
\str \hbox{Name} & \hfil\hbox{Defining polynomial} &\hfil \hbox{Intersection with $C_\vph$}\\
\hline\hline
\str k_{12} &\s_1^2\t_1^2 - \s_1\s_2\t_1\t_2 + \s_2^2\t_2^2
& D_{34}\,+\,D_{45}\,+\,3D_{14}\,+\,3D_{24}\\
\hline
\str m_{13} & \t_2
& D_{13}\,+\,D_{24}\\
\hline
\str k_{14} &  (\s_1\t_1 + \o^2\s_1\t_2 + \o\s_2\t_2)
(\s_1\t_1 + \o\s_1\t_2 + \o^2\s_2\t_2)
& D_{23}\,+\,D_{25}\,+\,3D_{12}+3D_{24}\\
\hline
\str m_{15} & \s_1
& D_{15}\,+\,D_{24}\\
\hline
\str m_{23} &  \t_1
&  D_{23}\,+\,D_{14}\\
\hline
\str k_{24} & (\s_1\t_1 + \o^2\s_2\t_1 + \o\s_2\t_2)
(\s_1\t_1 + \o\s_2\t_1 + \o^2\s_2\t_2)
& D_{13}\,+\,D_{15}\,+\,3D_{12}+3D_{14}\\
\hline
\str m_{25} & \s_2
&  D_{14}\,+\,D_{25}\\
\hline
\str m_{34} & \t_1-\t_2
& D_{12}\,+\,D_{34}\\
\hline
\str m_{35} & \s_1\t_1\,-\,\s_2\t_2
& D_{12}\,+\,D_{14}\,+\,D_{24}\,+\,D_{35}\\
\hline
\str m_{45} &  \s_1\,-\,\s_2
& D_{12}\,+\,D_{45}\\
\hline
\str l_1 &  \s_1^2 - \s_1\s_2 + \s_2^2
& D_{13}\,+\,D_{23}\,+\,D_{34}\,+\,D_{35}\\
\hline
\str l_2 & \t_1^2 - \t_1\t_2 + \t_2^2
& D_{15}\,+\,D_{25}\,+\,D_{35}\,+\,D_{45}\\
\hline
\end{tabular}
\capt{4.0in}{Divs}{The meromorphic functions on $C_\vph$ that arise as irreducible factors of the quantities
$\p_{ij}^5/p_c$ discussed in \sref{plmap}.}
\end{center}
\end{table}
\goodbreak
The other two polynomials are
$$
k_{24}\,:=\,\s_1^2\t_1^2 - \s_1\s_2\t_1^2 + \s_2^2\t_1^2 - \s_1\s_2\t_1\t_2 - \s_1\s_2\t_1^2 + \s_1^2\t_1^2
$$
which is obtained from $k_{14}$ by $(\s_1,\s_2)\leftrightarrow (\t_1,\t_2)$, i.e.\ by applying $(12)(35)$ in $\cS_5$, and
$$
k_{12}\,:=\,\s_1^2\t_1^2 - \s_1\s_2\t_1\t_2 + \s_2^2\t_2^2\,=\,
(\s_1\t_1+\o \s_2\t_2)(\s_1\t_1+\o^2 \s_2\t_2)~.
$$

A computation, similar to the one above, shows that
$$
(k_{12}=0)\,\cap\,C_\vph \,=\,D_{34}\,+\,D_{45}\,+\,3D_{14}+3D_{24}~.
$$

\tref{Divs} gives the zero divisors of these $12$ polynomials.
\begin{table}
\def\bigstr{\vrule height18pt depth12pt width0pt}
\def\str{\vrule height15pt depth10pt width0pt}
\begin{center}
\begin{tabular}[H]{|>{$}c<{$} |>{$~}l<{~$} || >{$}c<{$} | >{$~}l<{~$}|}
\hline
\multispan4{\bigstr\vrule\hfil\large  Factorization of the $p_{ij}:=\pi_{ij}^5/p_c$ \hfil\vrule}\\
\hline
\str \hbox{Name} & \hfil\hbox{Factorization} & \hbox{Name} & \hfil\hbox{Factorization} \\
\hline\hline
\str p_{12} &m_{34}\,m_{35}\,m_{45}\,k_{14}\,k_{24}\quad
& p_{24} &m_{13}\,m_{15}\,m_{35}\,k_{12}\,k_{14} \\
\hline
\str p_{13} & m_{13}^4\,m_{25}\,m_{45}\,k_{24}\,l_1
& p_{25} & m_{13}\,m_{34}\,m_{35}^4\,k_{14}\,l_2
\\
\hline
\str p_{14} & m_{23}\,m_{25}\,m_{35}\,k_{12}\,k_{24}
&p_{34} & m_{15}\,m_{25}\,m_{34}^4\,k_{12}\,l_1
\\
\hline
\str p_{15} & m_{15}^4\,m_{23}\,m_{34}\,k_{24}\, l_2
& p_{35} & m_{35}^4\,l_1\,l_2
\\
\hline
\str p_{23} &  m_{15}\,m_{23}^4\,m_{45}\,k_{14}\,l_1
&  p_{45} &  m_{13}\,m_{23}\,m_{45}^4\,k_{12}\,l_2
 \\
\hline
\end{tabular}
\capt{5.5in}{divpl}{Factorizations of the $p_{ij}$ in terms of the functions of the previous table.}
\end{center}
\end{table}
In \tref{divpl} we list the factorizations of the bidegree $(6,6)$ polynomials $p_{ij}=\pi_{ij}^5/p_c$ on $\IP^1{\times}\IP^1$. 
Using Table \ref{Divs} one can compute the divisors 
$(p_{ij}=0)\cap C^0_\vph$. 
For example, one finds that
\beq\begin{split}
\frac{\pi_{35}^5}{p_c}~&=~\frac{(\s\t)^4(1-\s)(1-\s+\s^2)(1-\t)(1-\t+\t^2)(\s\t-1)^5}
{\s^4\t^4(\s-1)(\t-1)(\s\t-1)}\\[10pt]
&=~(\s\t-1)^4(1-\s+\s^2)(1-\t+\t^2)~,
\end{split}\notag\eeq
so, after homogenizing and using the notation from \tref{Divs}, we get:
$$
p_{35}\,=\,m_{35}^4l_1l_2~.
$$
Thus we get: 
$$(p_{35}=0)\cap C_\vph\,=\,
4(D_{12}+D_{14}+D_{24}+D_{35})+
D_{13}+D_{15}+D_{23}+D_{25}+D_{34}+2D_{35}+D_{45}~.
$$
In this way one can determine the divisor $(p_{ij}=0)\cap C_\vph$ for all $ij$.
One finds, quite remarkably, that they can be written as:
$$
(p_{ij}=0)\cap C_\vph\,=\,D_b\,+\,5D_{ij}~,
$$
where the divisor $D_b$ does not depend on ${ij}$ and is given by:
$$
D_b\,=\,
4(D_{12}+D_{14}+D_{24})+
D_{13}+D_{15}+D_{23}+D_{25}+D_{34}+D_{35}+D_{45}~,
$$
so it is the sum of the $10$ divisors $D_{ij}$, but the ones corresponding to the singular points of $C^0_\vph$ have multiplicity four. 

Given $l\in \tC_\vph$, at least one of the $\pi_{ij}(l)$ is obviously non-zero. 
Thus the zeroes of the common factor $p_c$ as well as the contribution 
coming from the common zeroes of all $p_{ij}$'s are artefacts of the parametrization.
These common zeroes are the $10\cdot 2=20$ points in $D_b$ and these  correspond to van Geemen lines. To find the fifth powers of the Pl\"ucker coordinates of these points, one must take a limit on the curve $C_\vph$ (alternatively, one can use 
the explicit parametrizations of these lines given in equation (\ref{VanGone})).
The surface parametrizing the families of lines in all $\cM_\psi$ is thus a 
$(\IZ/5\IZ)^3$-covering of the blow up of $\dP_5$ in the 20 points~of~$D_b$.

An important consequence of our computations is that 
the meromorphic functions $\pi_{ij}^5/\pi_{pq}^5$ have zeroes of order 5 in the two points in $D_{ij}$ and have poles of order five in the two points in $D_{pq}$ since the (apparent) common zeroes of both cancel. These poles and zeroes correspond to van Geemen lines.

To be precise,
if $l$ is a line which has $\pi_{ij}(l)=0$, then it also has $\pi_{pq}(l)\neq0$ for some $pq$ and thus $\pi_{ij}^5(l)/\pi_{pq}^5(l)=0$, which shows that $l$ 
corresponds to a point in $D_{ij}$ on $C_\vph$, 
hence $l$ is a van Geemen line.

\subsection{The 125:1 cover $\tC_\vph$ of  $C_\vph$}\label{cover125}
We will now describe the Riemann surface $\tC_\vph$ more precisely.
We will use the construction of Riemann surfaces by means of polynomials 
$g_nT^n+g_{n-1}T^{n-1}+\ldots+g_0$ where the $g_i$ are meromorphic functions on a given Riemann surface. For example, the Fermat curve defined by $X^n+Y^n+Z^n=0$ in $\IP^2$ is the Riemann surface of the polynomial 
$T^n+(x^n+1)$ where $x$ is the meromorphic function on $\IP^1$ 
which gives the projective coordinate.

As we showed in section \ref{plmap},
the meromorphic function $\pi_{ij}^5/\pi_{pq}^5$, 
viewed as function on $C_\vph$,  
has zeroes of order 5 in the two points in $D_{ij}$ 
and it has poles of order five in the two points in $D_{pq}$ and is holomorphic, with no zeroes, on the rest of $C_\vph$.
We define the following meromorphic functions on $\tC_\vph$ and $C_\vph$ respectively:
$$
f_{ij}\,:=\,\pi_{ij}/\pi_{45},\qquad 
g_{ij}\,:=\,(\pi_{ij}/\pi_{45})^5\,=\,p_{ij}/p_{45}~.
$$
Notice that $f_{ij}/f_{pq}=\pi_{ij}/\pi_{pq}$, so we get all quotients of the Pl\"ucker coordinates from these $f_{ij}$. Obviously, $f_{ij}$ is a root of the polynomial 
$T^5-g_{ij}$. The other roots of the polynomial are the $\zeta^af_{ij}$ with $a=1,\ldots,4$.

The Riemann surface of this polynomial can be described as follows. Choose a coordinate  neighbourhood $U_x$ which  biholomorphic to a disc $\Delta\subset\IC$, with $0\in \Delta$ and local complex chart $z_x:U_x\rightarrow\Delta$ with $z_x(x)=0$.
If $x\in C_\vph$ and $g_{ij}$ has no poles or zeroes on $U_x$, this Riemann surface is locally the disjoint union of $5$ copies of $U_x$.
If $x$ is a zero of  $g_{ij}$, we can write 
$$
g_{ij}\,=\,z_x^5(1\,+\,a_1z_x\,+\,a_2z_x^2\,+\,\ldots )~.
$$ 
Restricting the open subset, we may assume that $1+a_1z_x+\ldots=h^5$ for a holomorphic function $h$ on $U_x$ without zeroes and poles.
On $U_x$ the polynomial is $T^5-(z_xh)^5=\prod_a(T-\zeta^az_xh)$, showing that the Riemann surface is still a disjoint union of $5$ copies of $U_x$. Another way to argue is that the subset $\{(z,t)\in\Delta^2:t^5=z^5\}$ is  a local model of the Riemann surface. This local model must be blown up in $(0,0)$ in order to get a smooth model.
For the poles of $g_{ij}$, which also have multiplicity five, one finds similarly that 
the Riemann surface is a disjoint union of $5$ copies of $U_x$. We refer to \cite{ForsterRiemannSurfaces} for these constructions of Riemann surfaces.

Thus the fact that each zero and pole of $g_{ij}$ has multiplicity five guarantees that the Riemann surface $\cX_{ij}$ of the polynomial $T^5-g_{ij}$ is an unramified covering of $C_\vph$. Since the $f_{ij}$ are meromorphic on $\tC_\vph$, there must exist holomorphic maps 
$$
\tC_\vph\,\longrightarrow\,\cX_{ij}\,\longrightarrow\,C_\vph
$$
with the first map of degree $25$. The second map, of degree $5$, is obtained from the polynomial 
$T^5-g_{ij}$. By \Mustata's results, $\tC_\vph$ is a connected Riemann surface, hence also $\cX_{ij}$ is connected. (Another way to see this is to notice that otherwise the polynomial 
$T^5-g_{ij}$ would be reducible. As its roots are the $\zeta^af_{ij}$, 
this would imply that there would be a meromorphic function $h_{ij}$ on $C_\vph$ with $h_{ij}^5=g_{ij}$. Then $h_{ij}$ would have poles, with multiplicity one, in only two points. Thus $C_\vph$ would be hyperelliptic. This is not the case, as the
map to $\IP^5$ induced by $\Phi$ is the canonical embedding of $C_\vph$.)

This construction can be iterated by considering the polynomial $T^5-g_{pq}$ on $\cX_{ij}$, for example, or by considering the fiber product of the Riemann surfaces $\cX_{ij}$ and
$\cX_{pq}$ over $C_\vph$. The main result is that $\tC_\vph$ can be obtained with this construction from three suitably choosen $g_{ij}$, for example, the $g_{i5}$, $i=1,2,3$. We have already remarked that the covering is unramified over points in the $D_{ij}$. Over each such point, we have found $125$ van Geemen lines.

Unramified covers with group $\cG\cong (\IZ/5\IZ)^3$ 
correspond to normal subgroups $K\subset \pi_1(C_\vph)$ 
of the fundamental group of $C_\vph$ with quotient 
$\pi_1(C_\vph)/K\cong \cG$. 
We will discuss an algebro-geometric approach to the covers with line bundles in \sref{restrictionmap}.
\newpage
\section{Special members of the Dwork pencil}\label{singularmanifolds}
\subsection{The case $\psi=0$, $\vph=\infty$}
Due to the relation defining $\vph^2$ in terms of $\psi$, the case $\psi=0$ corresponds to $\vph=\infty$. The quintic threefold $\cM_\psi$ is the Fermat quintic, and we already discussed the lines on this threefold in the introduction. The curve $C_\infty$ is the union of the $10$ lines on $\dP_5$. Now we would like to describe $\tC_\infty$ in more detail, using the description of the general $\tC_\vph$ as a 125:1 covering of $C_\vph$ defined by the polynomials $T^5-g_{ij}$, with $g_{ij}=p_{ij}/p_{45}$.

First of all, we restrict our attention to the line $E_{15}\subset C_\infty$, which corresponds to the curve $\sigma=0$ in $\IP^1{\times}\IP^1$. Putting $\s=0$ in all
$p_{ij}$, we notice first of all that as $m_{15}=\s$, we get:
$$
p_{15}\,=\,p_{23}\,=\,p_{24}\,=\,p_{34}\,=\,0\qquad\mbox{on}\quad E_{15}~.
$$
The restrictions of the other $p_{ij}$ are easy to compute, one finds:
$$
p_{12}\,=\,p_{25}\,=\,(\t-1)(\t^2-\t+1),\quad 
p_{23}\,=\,p_{35}\,=\,\t^2-\t+1,\quad
p_{45}\,=\,-p_{14}\,=\,(\t+1)(\t^2-\t+1)~.
$$
Thus we get a 25:1 covering of $E_{15}\cong\IP^1$, with coordinate $\t$, given by the two polynomials
$$
T^5\,+\,\t\,-\,1,\qquad U^5\,+\,\t\,+\,1~.
$$
The first polynomial has a zero, of order one, in $\t=1$ and a pole, of order one, in $\t=\infty$. The Riemann surface we get is the 5:1 cyclic cover of $\IP^1$ 
totally branched over these points. In particular it is a $\IP^1$ with coordinate $t$ satisfying $t^5=\tau-1$. Substituting this in the other polynomial, 
we get the polynomial
$U^5+t^5+2$, which defines (up to rescaling) the degree 5 Fermat curve. 
Using the $\cS_5$ action, we find that $\tC_\infty$ 
is the union of 10 Fermat curves of degree 5, as expected.
 
Now we consider the lines parametrized by the component of 
$\tC_\infty$ lying over $E_{15}$.
We already observed that for a line $l$ in this component we have $\pi_{ij}(l)=0$
for $ij=23,24,34$ since these $p_{ij}$ are zero.  So if $l$ is spanned by $x,y$ then 
$(x_2,x_3,x_4)$ and $(y_2,y_3,y_4)$ are linearly dependent, hence we may assume that $y=(y_1,0,0,0,y_5)$. As $l\subset\cM_0$ we get $y_1^5+y_5^5=0$, so we can put $y_1=1$, $y_5=\zeta^a$.
As also $p_{15}=0$, the vectors $(x_1,x_5)$ and $(y_1,y_5)$ are dependent and substracting a suitable multiple of $y$ from $x$ we see that $x=(0,x_2,x_3,x_4,0)$, and still $l=\langle x,y\rangle$. 
As $l\subset \cM_0$ we get $x_2^5+x_3^5+x_4^5=0$, which, after permuting the coordinates 2 and 5, gives the family of lines parametrized by the quintic Fermat curve described in the introduction. Thus we succeeded in recovering the lines on the Fermat quintic threefold with our description of $\tC_\vph$. 
\subsection{The cases $\psi^5=1$, $\vph^2=125/4$}
In case $\psi^5=1$, the threefold $\cM_\psi$ has 125 ordinary double points
and has been studied extensively in \cite{SchoenMumfordHorrocksBundle}.
For convenience, we will take $\psi=1$, $\vph=5\sqrt{5}/2$,
the other cases are similar.
A computation shows that the corresponding curves $C_\vph^0$ 
acquire 6 more ordinary double points. Since each double point lowers the genus by one, the desingularizations, which we denote by $\widehat{C}_\vph$ to distinguish them from the (singular) curves $C_\vph$ in $\dP_5$, are thus isomorphic to $\IP^1$. 

The $\pi_{ij}^5$, which are functions on $C_\vph^0$, are now fifth-powers
of functions on $C_\vph^0$.
Hence the 125:1 cover $\tC_\vph$ of $C_\vph$, given by the polynomials 
$T^5-(\p_{ij}/\pi_{45})^5$,
is the union of  125 copies of $\widehat{C}_\vph\cong \IP^1$.

This corresponds to the fact that $\cM_1$ contains 125 quadrics, each isomorphic to $\IP^1{\times}\IP^1$.  Each quadric has two families of lines, given by the $\{x\}\times\IP^1$ and $\IP^1\times\{x\}$ where $x$ runs over~$\IP^1$. Thus we get $2\cdot 125$ families of lines parametrized by $\IP^1$ in $\cM_1$.  
These correspond to the components of the coverings $\tC_\vph$ 
of the $C_\vph$ with $\vph^2=125/4$.

We will first discuss the lines on one of the quadrics, denoted by $Z$, in $\cM_1$.
We also give explicitly a (complicated) map from 
$\IP^1$ to $C_\vph^0$
which is a birational isomorphism.
One can then check that the fifth powers of the Pl\"ucker coordinates 
are now indeed fifth powers on $C_\vph^0$, hence the cover $\tC_\vph\rightarrow C_\vph$ becomes reducible.

The threefold $\cM_1$ has $125$ ordinary double points, 
they are the orbit of the point $q:=(1,1,1,1,1)$ under the action of $\cG$.
In the paper \cite{SchoenMumfordHorrocksBundle} it is shown that there are $125$ hyperplanes 
(i.e.\  linear subspaces $\IP^3\subset\IP^4$), which form one $\cG$-orbit, 
each of which  cuts $\cM_1$ in a smooth quadratic surface and a cubic surface. 
To see such a hyperplane, one writes the
equation (\ref{DworkPencil}) for $\cM_1$ as a polynomial in 
the elementary symmetric functions in $x_1,\ldots,x_5$:
$$
s_1\,:=\,\sum_{i=1}^5x_i~,\qquad s_2\,:=\,\sum_{i<j}x_ix_j~,\quad 
\ldots~,\quad s_5\,:=\,x_1x_2x_3x_4x_5~.
$$
The equation is then
$$
\cM_1\,:\qquad 
s_2s_3\,+\,s_1\left(s_4-s_2^2-s_1s_3+s_1^2s_2- \smallfrac{1}{5}s_1^4\right)\,=\,0~.
$$
Thus the hyperplane $H$ defined by $s_1=0$ cuts $\cM_1$ in the surface defined by
$s_2s_3=0$. One verifies that the quadric $Z\subset \cM_1$ defined by 
$s_1=s_2=0$ is a smooth quadric in $H\cong\IP^3$.
In $H$ we have $x_5=-(x_1+\ldots+x_4)$, hence $2s_2=(\sum x_i)^2-\sum x_i^2$
restricts to 
$-2\sum_{i\leq j} x_ix_j$.
Hence
$$
Z\,\cong\,\left\{(x_1,\ldots,x_4)\in\IP^3:\;\sum_{i\leq j} x_ix_j\,=\,0\,\right\}~.
$$
First of all, we are going to find the van Geemen lines in $Z$. 
Recall that these are the lines which are fixed under an element of order three 
in $\cS_5$. Taking this element to be $(123)$,
we thus try to find a constant $b$ such that the line, in $H$, parametrized by
$$
u\,\big(1,\o,\o^2,0,0 \big)\,+\,v\,\big(1,1,1,b,-(b+3)\big) ~=~
\big(u+v,\o u+ v,\o^2u+v,bv,-(b+3)v\big)~,
$$
lies in $Z\subset X$. Next we impose $s_2=0$ and we find the condition:
$$
b^2\,+\,3b\,+\,6\,=\,0,\qquad\mbox{hence}\quad
b_{\pm}\,=\,\frac{-3\pm\sqrt{-15}}{2}~,
$$
and we get two lines $l_\pm$ in $Z$ which meet in the point $(1,\o,\o^2,0,0)$.
From this one easily finds the other van Geemen lines on $Z$. 

The surface $Z$ is a nonsingular quadric in $\IP^3$ hence is isomorphic to $\IP^1{\times}\IP^1$. We wish to parametrize $Z$ and this parametrization is simplified by making an appropriate choice of coordinate on the first $\IP^1$, which we regard as the curve $C_\vph$ that parametrizes the lines. The group $\cA_5$, which is isomorphic to the icosahedral group, acts on $C_\vph$ and it is convenient to choose a coordinate $z$ adapted to this action. In the standard discussions of the automorphic functions of the icosahedral group~\cite{HigherTranscendentalFunctions, Forsyth}  one considers the projection of an icosahedron on to the circumscribing sphere and then the further projection of the image onto the equatorial plane, taking the south pole as the point of projection. There are thus two natural choices of coordinates, depending on whether the orientation of the icosahedon is chosen such that the south pole coincides with a vertex or the image of the center of a face. The standard treatments place a vertex at the south pole. We shall refer to this choice of coordinate, $w$, as the icosahedral coordinate. It can be checked that the 10 van Geemen lines of $C_\vph$ correspond to projection onto the circumscribing sphere of the centers of the faces of the icosahedron, or equivalently, to the vertices of the dual dodecahedron. For our purposes, it is therefore natural to work with a `dodecahedral coordinate' $z$ that corresponds to aligning the icosahedron such that the south pole of the circumscribing sphere corresponds to a vertex of the dual dodecahedron. The two coordinates may be chosen such that the relation between them is
$$
\o z~=~\frac{w_\infty\, w +1}{w\, -\, w_\infty}~,
$$
where $w_\infty$ denotes the $w$-coordinate of the dodecahedral vertex at the north pole of the circumscribing sphere. This can be chosen to be
$$
w_\infty~=~\frac14\left( 3+\sqrt{5} + \sqrt{6(5+\sqrt{5})}\right)~.
$$

It is convenient to fix a primitive $15$-th root of unity $\eta=e^{2\pi i/15}$. Then we also have a fifth and a third root of unity, $\z,\o$ respectively, and expressions for $\sqrt{5}$ and $w_\infty$:
$$
\zeta\,:=\,\eta^3,\qquad \o\,:=\,\eta^5,\qquad
\sqrt{5}\,=\,1+\zeta+\zeta^{-1},\qquad w_\infty\,=\, -2\eta^7 + \eta^5 - \eta^4 + \eta^3 - \eta + 2~. 
$$

Using the van Geemen lines, it is easy to find the following parametrization $\U:\IP{\times}\IP^1\rightarrow Z$,
where we reinstate the $x_5$ coordinate for symmetry reasons,
{\renewcommand{\arraystretch}{1.3}
$$
\U:\,(z,u)\,\longmapsto 
\left(\begin{array}{c}
x_1\\ x_2\\x_3\\x_4\\x_5
\end{array}\right)
\,=\,
\left(\begin{array}{r >{\hskip-7pt}l  >{\hskip-7pt}l  >{\hskip-7pt}l}
    -(b+3)\,c u z &&& +\,5b\\
   c u z& +\, 5 u &  +\, \o dz &+\, 5\\
   c u z&+\, 5 \o u &+\, d z& +\, 5\\
 b\,c u z&&&-\,5(b+3)\\
    c u z& +\,5\o^2 u& +\,\o^2 d z& +\, 5
\end{array}\right)~,
$$
}
where we make use of the following coefficients:
$$
\begin{array}{rcl}
b&:=&-\eta^7 +\eta^5-2 \eta^4+ \eta^3 -  \eta^2 -2 \eta~,\\[3pt]
c&:=&-2 \eta^7 + \eta^5 - 2 \eta^4 + 2 \eta^3 - 2 \eta^2 - 2 \eta + 2~,\\[3pt]
d&:=&-10 \eta^7 + 10 \eta^3 - 10 \eta^2 + 5~,
\end{array}
$$
in particular, $b^2+3b+6=0$.
For fixed $z$, we have a map $\IP^1\rightarrow Z$ 
whose image is a line $l_z$ in $Z$ parametrized by $u$. 
One can check that the action of $\cA_5$, which has generators of order 2, 3 and 5, on the coordinates
$x_1,\ldots,x_5$ corresponds to the action of the following 
M\"obius transformations:
$$
M_2(z)\,:=\,-1/z,\qquad 
M_3(z)\,:=\,\o z,\qquad 
M_5(z)\,:=\,
\frac{(\z w_\infty^2+1)\, z + (\z-1)\, \o^2 w_\infty}{(\z-1)\, \o\, w_\infty\,z + (\z+w_\infty^2)}~,
$$
where the order 5 transformation $M_5$ is simply the transformation $w\to \z w$, when written in terms of the icosahedral coordinate. 

The polynomial whose roots, together with $z=0$ and $z=\infty$, correspond to the dodecahedral vertices is
$$
8\, z^{18} - 57 \sqrt{5}\, z^{15} - 228\, z^{12} - 494 \sqrt{5}\, z^9 + 228\, z^6 - 57 \sqrt{5}\, z^3 - 8~.
$$
The van Geemen lines correspond to the nine pairs of roots $\{z_{*},-1/z_{*}\}$ together with $\{0,\infty\}$.

For the M\"obius transformation $M_k$ one has 
$$
l_{M_k(z)}\,=\,g_k(l_z)~,\qquad\mbox{with}\quad
l_z\,:=\,\{\U(z,u)\,\in\,\IP^4:\,u\in\,\IP^1\}~,
$$
where, in this context,
$$
g_k\,=\,(14)(25),\;(253),\;(54321)~~~\text{for}~~~k\,=\,2,3,5~.
$$
The orbit of the line $l_z$, with $z=0$, which is a van Geemen line fixed by 
$(253)$, consists of $20$ van Geemen lines.

On $\cM_1$ there are also lines fixed by an element of order five in $\cA_5$. 
These are the lines that cause the extra double points on $C_\vph$. 
The element of order five $(12345)\in\cA_5$ has five isolated fixed points in $\IP^4$, four of which lie on  $Z=H\cap \cM_1$, in fact they are singular points of $\cM_1$.
They are, for $i=1,\ldots,4$:
$$
q_i\,:=\,(\, \zeta_i^j \,)_{1\leq j\leq 5}~,\qquad
\big(\{q_1,q_2,q_3,q_4\}\subset \mbox{Sing}(\cM_1)\big)~.
$$
One easily checks that $\lambda q_1+\mu q_2+\nu q_3$ lies on $Z$
only if $\mu=0$ or $\nu=0$ and thus the lines $\lambda q_1+\mu q_2$ and 
$\lambda q_1+\nu q_3$ do lie in $Z$. 
So the intersection of the $\IP^2$ spanned by $q_1,q_2,q_3$ with the quadric $Z$
consists of two lines, each of which is spanned by two nodes:
$$
\langle q_1,\,q_2,\,q_3\,\rangle\,\cap\,Z\,=\,
\langle q_1,\,q_2\rangle\,\cup\,\langle q_1,\,q_3\rangle~.
$$
Both of these lines are invariant under the $5$-cycle $(12345)\in\cS_5$, 
and similarly we get lines $\langle q_2,\,q_4\rangle$ and $\langle q_3,\,q_4\rangle$.
The two lines $\langle q_1,\,q_2\rangle$, $\langle q_3,\,q_4\rangle$ on $Z$
are fixed by the 5-cycle and they do not intersect, 
hence they are from the same ruling.
Applying $\cA_5$, we get 12 lines, actually six pairs, with a stabilizer of order five
in each ruling. These create the 6 double points in~$C_\vph$.

We now briefly discuss the curve $C_\vph^0$ and a parametrization.
The curve $C_{\vph}^0$ has 6 more ordinary double points where $(\s,\t)$ 
take the values
$$
\left( \pm\smallfrac12 (1+\sqrt{5}),\, \pm \smallfrac12 (1+\sqrt{5}) \right)\, ,~
\left( \pm \smallfrac12 (1+\sqrt{5}),\,  \smallfrac12 (3-\sqrt{5})  \right)~~\text{and}~~
\left(\smallfrac12 (3-\sqrt{5}) ,\,  \pm \smallfrac12 (1+\sqrt{5})\right)\, ,
$$
where, in the first expression, the same sign is chosen for each component.

A parametrization of $C_\vph^0$ is given by
$$
\IP^1\,\longrightarrow\,C_\vph^0,\qquad
z\,\longrightarrow\,(\s,\t)\,=\,\big(R_1(z),\,R_2(z)\big)
$$
with rational functions
\beq\begin{split}
R_1(z,\vph)~&=~\frac{20z^4 + (15 - 2\vph) z^3 + 3(5 + 2\vph)z^2 - (15 - 2\vph)z + 20}
{6z(5z^2 + 2\vph z - 5)}~,\\[10pt]
R_2(z,\vph)~&=~\frac{1}{R_1(-\o^2 z,-\vph)}~.
\end{split}\notag\eeq

In particular, one has $F_+\big(R_1(z),R_2(z)\big)=0$ for all $z$.

The $\cA_5$-action on $C^0_\vph$ lifts to an action of $\cA_5$ by M\"obius transformations on $\IP^1$.
Generating transformations are the same $m_k$ as given earlier 
and one has 
$$
\Big( R_1\big(m_k(z)\big),\,R_2\big(m_k(z)\big)\Big)\,=\,g_k\big(R_1(z),\,R_2(z)\big)~,
$$
where again,
$$
g_k\,=\,(14)(25),\;(253),\;(54321)~~~\text{for}~~~k\,=\,2,3,5~.
$$
This parametrization can be found using the $\cA_5$-action.
The coordinate function $\sigma$ is known to be the quotient map by the subgroup generated by $(12)(34)$ and $(13)(24)$. 
Since we require the map to be equivariant, the points in $z\in\IP^1$ which are zeroes/poles of $R_1$ must correspond to van Geemen lines and these are fixed points of order three. As also the fiber over $R_i^{-1}(1)$ must consist of such fixed points, the $R_i$ are easily found.

The points $z=0$, $z=\infty$ in $\IP^1$ are both mapped to the singular point
$(\infty,0)\in C^0_\vph\subset\IP^1{\times}\IP^1$. Using Table \ref{ExcCurves}
we find that these two points correspond to the divisor $D_{14}$. 
The M\"obius transformations $m_2$ and $m_3$ on $\IP^1$ fix the set $\{0,\infty\}$
and this is indeed consistent with the fact that the permutations $(14)(25)$ and $(253)$ fix the index set $\{1,4\}$ and thus fix the divisor $D_{14}$ on $C_\vph$.
As $\cA_5$ acts transitively on the set of $20$ points which are in $\cup D_{ij}$,
the $\cA_5$-orbit of $0\in\IP^1$ consists of the $20$ points which map bijectively to this set.

The two fixed points of $m_5$ in $\IP^1$ map to the singular point 
$\big(\!-\smallfrac12 (1+\sqrt{5}),\, \smallfrac12 (3-\sqrt{5})\big)$ in~$C^0_\vph$.
The $\cA_5$-orbit of any of these points is a set of 12 points, each a fixed point of an order 5 element in $\cA_5$, which maps to one of the other six `extra' singular points of $C_\vph^0$.

Finally we consider the fifth powers of Pl\"ucker coordinates of the lines 
parametrized by $C_\vph$. These functions are given by the polynomials
$p_{ij}(\s_1,\s_2,\t_1,\t_2)$ listed in Table \ref{divpl}.
We pull them back to $\IP^1$ along the parametrization, so we 
take homogeneous coordinates $(z,w)$ on $\IP^1$ and we
consider the $10$ polynomials
$$
\tilde{p}_{ij}(z,w)\,:=\,
p_{ij}(R_{1,n}(z,w),\,R_{2,n}(z,w),\,R_{1,d}(z,w),\,R_{2,d}(z,w))
$$
where we introduced
$$
R_i(z,w)\,:=\,w^4R_i(z/w)\,=\,R_{i,n}(z,w)/R_{i,d}(z,w)
$$
and the $R_{i,n}(z,w)$, $R_{i,d}(z,w)$ are 
homogeneous polynomials of degree four.
As $p_{ij}$ is homogeneous of bidegree $(6,6)$,
the polynomials $\tilde{p}_{ij}$ are homogeneous of degree 
$6\cdot 4+6\cdot 4=48$.
Using the results from Section \ref{plmap} and the fact that we are now on a $\IP^1$,
the divisor $D_b$, which has degree $38$, is now defined by a homogeneous polynomial $\tilde{p}_b(z,w)$ of degree $38$. 
Each of the $\tilde{p}_{ij}$ is divisible by $\tilde{p}_b$ with quotient $\tilde{q}_{ij}$ which is homogeneous of degree $10$. 
We know that each $\tilde{q}_{ij}$ has two zeroes, with multiplicity $5$, 
in the two points of $D_{ij}$. We checked that this is indeed the case.

In fact, we verified that the point $(\ldots,q_{ij}(z),\ldots)\in \IP^9$ is the point 
$(\ldots,\pi_{ij}(l_z)^5,\ldots)$, the point whose coordinates are the fifth powers of the Pl\"ucker coordinates of the line $l_z\subset \cM_1$. Each $\pi_{ij}(l_z)$ is easily seen to be a quadratic polynomial in $z$.

Thus the parametrizations of $Z$ and
$C^0_\vph$ are compatible and the $2\cdot 125$ families of
lines on the 125 quadrics in $\cM_1$ 
are the limits of the two families of lines on the general $\cM_\psi$.

In particular, the curve $\tC_\vph$ is now reducible, having $125$ components, the desingularization of each component is a $\IP^1$.
Moreover, the Pl\"ucker map from $\tC_\vph$ to the Grassmannian in $\IP^9$ 
is given by $10$ degree two polynomials on each of the components, the fifth power of these polynomials are the $\tilde{q}_{ij}$.

We checked that the $125$ components of $\tC_\vph$ intersect as follows. On each component there are the 12 fixed points of certain elements of order five in $\cS_5{\rtimes}\cG$. In each such point, exactly two components meet and 
moreover, distinct components only meet in such fixed points. 
Thus $\tC_\vph$ has $(125\cdot 12)/2=750$ ordinary double points.

This allows us to compute the arithmetic genus $p_a(\tC_\vph)$ of the curve $\tC_\vph$, since
$1-p_a(\tC_\vph)$ is the Euler characteristic of the intersection graph of the components of $C_\vph$. This graph has $125$ vertices and 
$750$ edges, hence it has Euler characteristic $125-750=-625$. Thus indeed $p_a(\tC_\vph)=626$, as expected.

\subsection{The case $\psi=\infty$, $\vph^2=-3/4$}\label{psiinfty}
In case $\psi=\infty$ the threefold $\cM_\psi$ is defined by $x_1\cdots x_5=0$, 
so it is the union of $5$ hyperplanes.
The corresponding curves $C_\vph^0$ become reducible, 
in fact the polynomial $F_+$ has 5 factors for these values of $\vph$:
$$
F_+\,=\,(\s +\o^2)(\t+\o^2)(\s\t+\o)(\s\t + \o\s + \o^2)(\s\t +\o\t + \o^2)~,
$$
and $F_-$ is obtained by $\o\leftrightarrow\o^2$. 
These factors are also factors of $l_1$, $l_2$, $k_{12}$, $k_{14}$ and $k_{24}$ respectively (see Table \ref{Divs}). 
The components of these $C_\vph$, and their classes in $\Pic(\dP_5)$, 
are discussed at the end of Section \ref{picdp5}. 
Each component of $C_\vph$ parametrizes lines in one of the hyperplane $x_i=0$ in $\cM_\psi$, these $x_i$ are $x_3$, $x_5$, $x_4$, $x_2$ and $x_1$ respectively.
The cover $\tC_\vph\rightarrow C_\vph$ is non-trivial in this case and we will not analyze it any further here.

For example, assume that we are on the component where $\s=-\o^2$.
Recall that the $p_{ij}(-\o^2,\t)$ are, upto a common factor, the 
$\pi_{ij}^5(-\o^2,\t)$. These polynomials are listed in Table \ref{divpl} and
one finds that 
$$
p_{ij}(-\o^2,\t)\,=\,0\qquad\mbox{for}\quad ij\in\{13,23,34,35\}~.
$$ 
Thus this component of $C_\vph$ parametrizes lines $l$ 
which have  $\pi_{i3}(l)=0$ for all $i$. 
Such a line $l$ lies in the hyperplane $x_3=0$, 
because else we may assume that $l=\langle x,y\rangle$ with $x_3\neq 0$, 
in which case we may assume that $y_3=0$ and moreover one $y_j$, $j\neq 0$ must also be non-zero, but then $\pi_{j3}\neq 0$.

Moreover, after dividing the six non-zero polynomials $p_{ij}(-\o^2,\t)$ 
by a common factor of degree $3$, 
the quotients $q_{ij}(\t)$ are degree two polynomials in $\t$. 
Define
$$
n_1\,:=\,\t+(\o-1)/3,\quad n_2\,:=\,\t+\o-1,\quad
n_4\,:=\,\t+\o+1,\quad n_5\,:=\,\t-\o-1~.
$$
Then we have, for $i,j\in \{1,2,4,5\}$ and $q_{ij}=0$ else, that for certain $c_{ij}\in\IC$,
$$
q_{ij}\,=\,c_{ij}n_in_j,\qquad\mbox{and}\quad
(\ldots,p_{ij}(-\o^2,\t),\ldots)\,=\,(\ldots,q_{ij}(\t),\ldots)\quad(\subset\IP^9)~.
$$

\newpage
\section{Appendix}
\subsection{The Picard group of $\dP_5$}\label{picdp5}
We recall the basic facts on the geometry of the quintic del Pezzo surface $\dP_5$.
We will use some more advanced algebraic geometry in this section to put the results we found in a more general perspective.

It is most convenient to view $\dP_5$ as the blow up of $\IP^2$ in four distinct points, no three on a line.
One can then choose coordinates such that these four points are
$$
p_1=(1,0,0),\quad p_2=(0,1,0), \quad
p_3:=(0,0,1),\quad p_4:=(1,1,1)~.
$$
The blow up map $\IP^2\dashrightarrow \dP_5$ is given by the cubic polynomials which are zero in the $p_i$.

There is an obvious action of $\cS_3$ by automorphisms of $\dP_5$ induced by permutation of the coordinates.
The action of $\cS_3$ extends to a linear action of $\cS_4$, 
the subgroup of $\text{PGL}(3,\IC)=\Aut(\IP^2)$ which permutes the four points.
In fact, the map
$$
\sigma_{34}:\,\IP^2\,\longrightarrow\,\IP^2,\qquad (x:y:z)\,\longmapsto\,
(x-z,y-z,-z)
$$
fixes the first two points and exchanges the last two.
Finally the standard (birational) Cremona transformation
$$
\sigma_{45}:\,\IP^2\,\longrightarrow\,\IP^2,
\qquad (x,y,z)\,\longmapsto\,(x^{-1},y^{-1},z^{-1})\,=\,(yz,\,xz,\,xy)
$$
induces another automorphism of $\dP_5$, which together with the 
$\cS_4$ generates a group isomorphic to $\cS_5$ and $\cS_5=\Aut(\dP_5)$.

The quintic Del Pezzo surface $\dP_5$ has $10$ exceptional divisors, 
which we denote by $E_{ij}=E_{ji}$ with $1\leq i<j\leq 5$. The divisors $E_{i5}$ are the exceptional divisors over the points $p_i$ and the $E_{ij}$,  with $1\leq i<j\leq 4$ are (somewhat perversely, but this helps in understanding 
the intersection numbers) is the strict transform of the line 
$l_{ij}$ spanned by $p_k$ and $p_l$, with $\{i,j,k,l\}=\{1,2,3,4\}$.
So the pull-back of the line $l_{ij}$ in $\IP^2$ to $\dP_5$ has divisor $E_{ij}+E_{k5}+E_{l5}$ and, for example, $l_{12}$ is defined by $x-y=0$,
$l_{24}$ is defined by $y=0$. In particular we have
$$
E_{ij}\,=\,l\,-\,E_{k5}\,-\,E_{l5}\qquad(\in\,\Pic(\dP_5))~.
$$

With these conventions, the intersection numbers are
$$
E_{ij}^2\,=\,0~,\qquad E_{ij}E_{ik}\,=0\quad
\mbox{if}\quad\sharp\{i,j,k\}=3~,\qquad
E_{ij}E_{kl}\,=1\quad\mbox{if}\quad\sharp\{i,j,k,l\}=4~.
$$
The intersection graph of the $E_{ij}$ has $10$ vertices and $15$ edges, each vertex is on three edges. This graph is known as the Petersen graph and is presented in \fref{petersengraph}.
\begin{figure}[H]
\begin{center}
\includegraphics[width=2.5in]{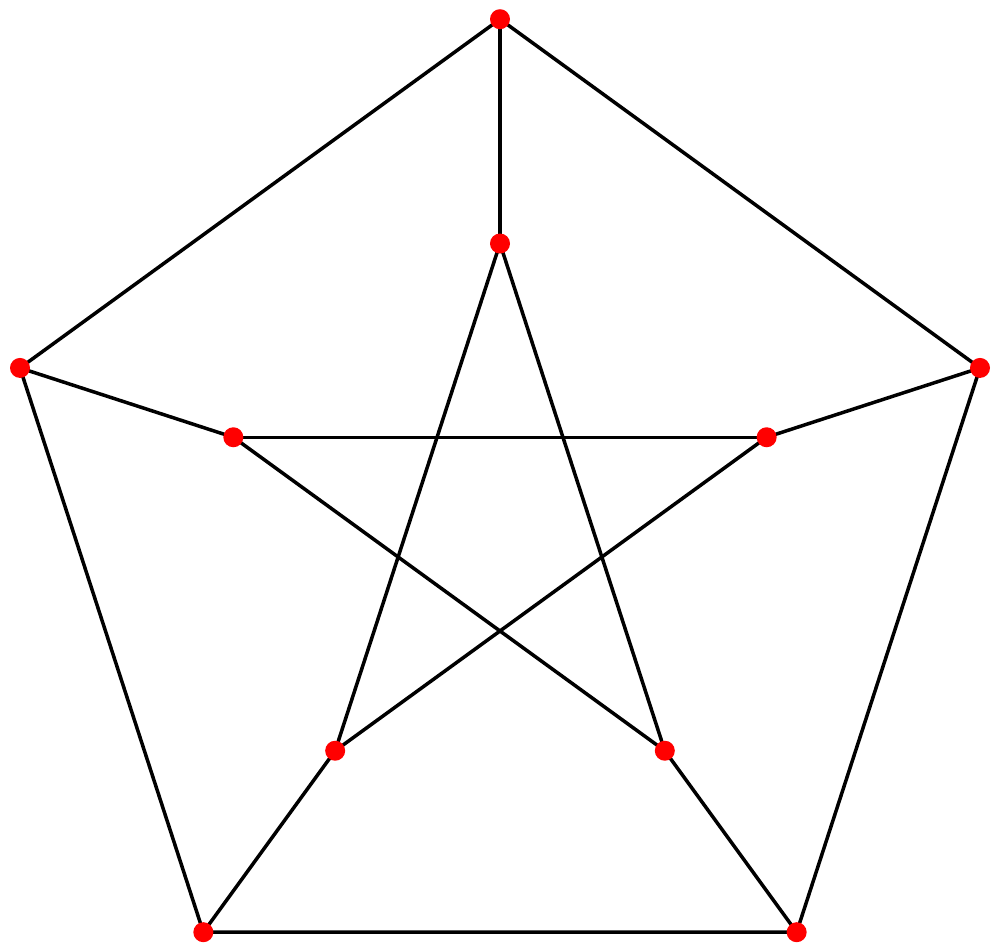}
\vskip10pt
\capt{5.5in}{petersengraph}{The Petersen graph, which summarizes the combinatorics of the intersections of the exceptional divisors $E_{ij}$. In the figure, the exceptional divisors correspond to vertices and their intersections correspond to the edges.}
\end{center}
\end{figure}
Let $l\in \Pic(\dP_5)$ be the class of the pull-back of a line in $\IP^2$.
One has $l^2=+1$.
Then the canonical class $K_{\dP_5}$ of $\dP_5$ is determined by 
$$
-K_{\dP_5}\,=\,3l\,-\,E_{15}\,-\,E_{25}\,-\,E_{35}\,-\,E_{45}
\qquad (\in\,\Pic(\dP_5))~,
$$
we have $(-K_{\dP_5})^2=9-4\cdot 1=5$.
In particular, the anti-canonical map of $\dP_5$ is induced by the cubics on the four nodes of $C_\vph$. One also has
$$
\Pic(\dP_5)\,=\, \IZ l\,\oplus\,\IZ E_{15}\,\oplus\,\IZ E_{25}\,\oplus\IZ E_{35}\,\oplus\,\IZ E_{45}~.
$$

The action of $\cS_5$ on $\Pic(\dP_5)$ is as follows. 
The permutations which fix $5$  are induced by linear maps on $\IP^2$
and thus act by fixing $l$ and permuting the indices of the $E_{ij}$.
The transposition $(45)$ is induced by the Cremona transformation. The pull-back of a line is a conic on $p_1,p_2,p_3$, thus
$\s_{45}^*l=2l-E_{15}-E_{25}-E_{35}$ and the image of the line on $p_i,p_j$ is the point $p_k$, with $\{i,j,k\}=\{1,2,3\}$ so 
$\sigma_{45}^*E_{k5}=l-E_i-E_j$. The point $p_4=(1,1,1)$ is mapped to itself, so $\sigma^*E_{45}=E_{45}$.

The Picard group of $\IP^1{\times}\IP^1$ is generated by the classes of the divisors $\s=0$ and $\t=0$. The holomorphic map 
$\Phi:\dP_5\rightarrow \IP^1{\times}\IP^1$ induces the pull-back homomorphism
$$
\Phi^*:\,\Pic(\IP^1{\times}\IP^1)\,\longrightarrow \, \Pic(\dP_5)~,
\qquad
\left\{\begin{array}{rcl}
(\s=0)&\longmapsto&E_{15}\,+\,E_{24}~,\\
(\t=0)&\longmapsto&E_{23}\,+\,E_{14}~.
\end{array}\right.
$$
Notice that the curve $E_{15}$ maps to $\s=0$ under $\Phi$, 
but the point $(0,\infty)$ on $\s=0$ is blown up, so its exceptional divisor $E_{24}$ contributes to $\Phi^*(s=0)$. 
In the standard basis of $\Pic(\dP_5)$ we have
$$
\begin{array}{rcl}
(\s=0)&\longmapsto&E_{15}\,+\,E_{24}\,=\,l\,-\,E_{35}~,\\
(\t=0)&\longmapsto&E_{23}\,+\,E_{14}\,=
\,2l\,-\,E_{15}\,-\,E_{25}\,-\,E_{35}\,-\,E_{45}~,
\end{array} 
$$
showing that $\s,\t$ are related to lines through the point $p_3$ and conics on the four points $p_1,\ldots,p_4$ in $\IP^2$.

The curves of bidegree $(n,m)$ in $\IP^1{\times}\IP^1$ 
pull-back along $\Phi^*$ to curves with class 
$n(l-E_{35})+m(2l-\sum_iE_{i5})$ in $\Pic(\dP_5)$.
In case a point $p$ which gets blown up is a point with multiplicity $r$
on a curve, its pull-back to $\dP_5$ is reducible. One component is the exceptional divisor over $p$, which has multiplicity $r$, and the other is called the strict transform of the curve. In particular, a curve of type $(2,2)$ which passes through all three points which get blown up with multiplicity one  thus has four components, and its strict transform has class
$$
2(l-E_{35})\,+\,2(2l-\sum_iE_{i5})\,-\,E_{12}\,-\,E_{14}\,-\,E_{24}~.
$$
Using that $E_{ij}=l-E_{k5}-E_{l5}$ we find this class is equal to:
$$
3l\,-\,E_{15}\,-\,E_{25}\,-\,E_{35}\,-\,E_{45}\,=\,-K_{\dP_5}~.
$$
Thus we see that the rational map $\Psi$ from Section \ref{dp5}
induces the anti-canonical map on $\dP_5$ and it is known that this map embeds $\dP_5$ in $\IP^5$.  

The curves $C_\vph^0$ have bidegree $(4,4)$ and they have multiplicity $2$ 
in each  of the three points which get blown up, hence their strict transform 
$C_\vph$ in $\dP_5$ has class $-2K_{\dP_5}$. 
It is easy to verify that the sum of the 10 exceptional divisors also has this class:
$$
\sum_{i<j}E_{ij}\,=\,6l\,-\,2E_{15}\,-\,2E_{25}\,-\,2E_{35}\,-\,2E_{45}\,=\,
-2K_{\dP_5}~.
$$
As the (reducible, singular) curve $\cup E_{ij}$ is also $\cA_5$-invariant, it should be in the Wiman pencil, and in fact it is $C_\infty$.

As each $C_\vph$ has class $-2K_{\dP_5}$, the intersection number of two such curves
is $(-2K_{\dP_5})^2=4\cdot 5=20$. Thus the Wiman pencil has 20 basepoints. We already
found $2\cdot 10=20$ points in the divisors $D_{ij}$ in the base locus, 
so the base locus is the union of the ten $D_{ij}$.

\subsection{The case $\psi=\infty$, $\vph^2=-3/4$ again}
We now give a more intrinsic description of the curves $C_\vph$ with $\vph^2=-3/4$ from Section \ref{psiinfty}. 

In Section  \ref{psiinfty}
we found that the curves $C_\vph$ with $\vph^2=-3/4$
have $5$ irreducible components, 
the first two of bidegree $(1,0)$, $(0,1)$ respectively, the last three of bidegree $(1,1)$. 
Only the last three curves pass through the base points $(1,1)$, $(0,\infty)$ and $(\infty,0)$: each contains two base points and each base point is on exactly two components. 
The class of the component of bidegree $(1,1)$ passing through 
$(1,1)$ and $(0,\infty)$ is
$$
(l\,-\,E_{35})\,+\,(2l\,-\,E_{15}\,-\,E_{25}\,-\,E_{35}\,-\,E_{45})\,-E_{12}\,-\,E_{24}\,=\,
l\,-\,E_{25}~,
$$
and similarly, the components passing through 
$(\infty,0)$, $(0,\infty)$ and $(1,1)$, $( \infty,0)$ are $l-E_{45}$ 
and $l-E_{15}$ respectively.
Thus the classes
of the strict transforms of these components are:
$$
l\,-\,E_{15},\quad l\,-\,E_{25},\,\quad\,l\,-\,E_{35},\quad l\,-\,E_{45},\,\quad
2l\,-\,(E_{15}+\ldots+E_{45})~.
$$
These five classes are in one orbit under $\cA_5$ (in fact under $\cS_5$), and they correspond to the 5 coordinates on $\IP^4$, as we also found in 
Section \ref{psiinfty}.

The $5$ components are not uniquely determined by their classes. 
In fact, each class determines a pencil of curves. 
The first four correspond to the pencil of lines through the point $p_i$
and the last to the pencil of conics on $p_1,\ldots,p_4$.
We will now use the action of $\cS_5$ on $\dP_5$ to find two specific curves in these pencils.

Each pencil is fixed by the subgroup, isomorphic to  $\cS_4$,
in $\cS_5$ which fixes the class. Thus an element $g$ of the $\cS_4$ corresponding to the
pencil maps a curve in the pencil to another curve in the pencil. For example,
the pencil of conics on $p_1,\ldots,p_4$ is fixed by the standard $\cS_4\subset \cS_5$.
A conic $C_{(\lambda,\mu)}$, with $(x,y)\in\IP^1$,  from this pencil has equation 
$$
\lambda (x-z)y\,+\,\mu(y-x)z\,=\,0~.
$$
There are three reducible conics in the pencil, $(x-z)y=0$, $(y-x)z=0$ and 
$(z-y)x=-(x-z)y-(y-x)z$, these consist of pairs of exceptional curves.
If $g\in \cS_4$ lies in the Klein subgroup $\cK:=\langle (12)(34),(13)(24)\rangle$, then
one verifies easily that $g$ maps $C_{(\lambda,\mu)}$ into itself.
Thus the action of $\cS_4$ on $\IP^1$ factors over quotient $\cS_4/\cK\cong\cS_3$.
The element $(123)\in \cS_5$ maps onto a generator of the subgroup $\cA_3$ of the quotient group $\cS_3$. It acts on the pencil as
$$
\lambda (x-z)y+\mu(y-x)z\,\longmapsto\,
\lambda (y-x)z+\mu(z-y)x\,=\,-\mu(x-z)y+(\lambda-\mu) (y-x)z,
$$
in particular, it has two fixed points $(\lambda,\mu)=(1,-\o), (1,-\o^2)$.
 
Thus the 5 classes above give $2\cdot 5=10$
curves in $\dP_5$ which we denote by $D_{ia}$, $D_{ib}$, $i=1,\ldots,5$. 
Upto permutations of $a,b$, the two curves  $\cup_i C_{ia}$ and $\cup_i C_{ib}$
are invariant under the action of $\cA_5$ and they are the
$C_\vph$, for $\vph^2=-3/4$.

It is interesting to notice that the action of $S_4$ on the five pencils shows that the
five maps from $C_\vph$ to $\IP^1$ they define are actually 
$(\IZ/2\IZ)^2$-quotient maps. For example,
from Table \ref{S5transfs}
one finds that $(12)(45), (14)(25)\in\cS_5$ act as
$$
(\s,\t)\,\longmapsto\, \left(\,\s,\,\frac{1}{\s\t}\right),\qquad
(\s,\t)\,\longmapsto\, \left(\,\s,\,\frac{\s\t-1}{\s(\t-1)}\right)
$$ 
on $\IP^1{\times}\IP^1$. Thus they 
fix $\s$, hence they act on the fibers of the projection map $(\s,\t)\mapsto \s$. 
So this projection map is invariant under the Klein subgroup $\langle (12)(34),(13)(24)\rangle$ of $\cS_5$. 
As the map has degree four, it follows that
the quotient of $C_\vph^0$ by this Klein subgroup is $\IP^1$, 
with quotient map $\s$.

\subsection{From $\IP^2$ to $\IP^1\times\IP^1$ and back}\label{P2blowup}
In Section \ref{dp5} we obtained $\dP_5$ as the blow up of $\IP^1{\times}\IP^1$
in three points. Blowing down the four exceptional curves $E_{15}$, $\ldots$, $E_{45}$ on $\dP_5$, we get $\IP^2$. The composition of these maps is a 
birational map between $\IP^1{\times}\IP^1$ and $\IP^2$.
To find it, we observe that $\Phi^*$ acts on the following divisors as:
$$
\begin{array}{cclcccl}
(\s=0)&\longmapsto&E_{15}\,+\,E_{24}~,&\qquad\qquad&
(\t=0)&\longmapsto&E_{23}\,+\,E_{14}~,\\
(\s=\infty)&\longmapsto&E_{14}\,+\,E_{25}~,&&
(\t=\infty)&\longmapsto&E_{13}\,+\,E_{24}~,\\
(\s=1)&\longmapsto&E_{12}\,+\,E_{45}~,&&
(\t=1)&\longmapsto&E_{12}\,+\,E_{34}~.
\end{array}
$$
Thus the function $\s$ on $\IP^1{\times}\IP^1$
corresponds to the pencil of lines in $\IP^2$ passing through the point
$p_3=(0,0,1)$ and in fact the meromorphic function $y/x$ on $\IP^2$ 
gives the same divisors on $\dP_5$.
Similarly $\t$ corresponds to the pencil of conics in $\IP^2$ 
passing through all four $p_i$ and its divisors match those of 
the meromorphic function $x(y-z)/y(x-z)$ on $\IP^2$.
Therefore the birational map from $\IP^2$ to $\IP^1{\times}\IP^1$ 
is given by
$$
\s\,=\,\frac{y}{x}~,\qquad \t\,=\,\frac{x(y-z)}{y(x-z)}~.
$$
As then $y=\s x$, one finds upon substitution in $\t=x(y-z)/y(x-z)$ and some manipulations that
$$
x:=\,\s\t - 1~,\qquad y\,:=\,\s(\s\t - 1)~,\qquad z\,:=\,\s(\t - 1)~.
$$
gives the inverse birational map.
These three polynomials are linear combinations of the polynomials $z_0,\ldots,z_5$ from Section \ref{dp5}, thus this map factors indeed over $\dP_5$.

It is amusing to verify that this works as advertised: 
take for example the curve defined by $\s=0$ on $\IP^1{\times}\IP^1$,
it maps to the exceptional divisor $E_{15}$ in $\dP_5$ according to Table 
\ref{ExcCurves}, and thus it should map to the point $p_1=(1,0,0)\in\IP^2$, 
which it does: $(x,y,z)=(-1,0,0)=(1,0,0)$.
Conversely, the line $l_{24}$ spanned by $p_1,p_3$ maps to $E_{24}$ in $\dP_5$
and next $E_{24}$ is mapped, according to the same table,  
to the point 
$(0,\infty)$ in $\IP^1{\times}\IP^1$. 
Indeed, $l_{24}$ is parametrized by $(a:0:b)$ and
thus $\s=y/x=0$ and $\t=x(y-z)/y(x-z)=-ab/0=\infty$.

The curve $C_\vph^0$ in $\IP^1{\times}\IP^1$ is defined by $F_+=0$.
We found polynomials $f_e,f_o$ in $x,y,z$ such that
$$
(xy(x-z))^4F_+\left(\frac{y}{x},\frac{x(y-z)}{y(x-z)}\right)\,=\,
\big(xy(x-y)\big)^2\big(f_e(x,y,z)\,-\,\vph f_o(x,y,z)\big)~.
$$
Thus the equation for the image of $C^0_\vph$ in $\IP^2$ is:
$$
f_e(x,y,z)\,-\,\vph f_o(x,y,z)\,=\,0~.
$$
This equation is homogeneous of degree six, it has an even and an odd part (under the action of $S_3$ which permutes the variables), where
$$
f_e\,=\,2s_1^2s_2^2\,-\,6s_1^3s_3\,-\,6s_2^3\,+\,19s_1s_2s_3\,-\,9s_3^2~,
$$
and the $s_i$ are the elementary symmetric function in $x,y,z$:
$$
s_1\,:=\,x\,+\,y\,+\,z~,\qquad
s_2\,:=\,xy\,+\,xz\,+\,yz~,\qquad
s_3\,:=\,xyz~.
$$
The odd part is
$$
f_o\,:=\,2xyz(x-y)(x-z)(y-z)~.
$$
In particular, any odd element in $\cS_5$ maps $C_\vph$ to $C_{-\vph}$, 
as we have already seen. The singular locus of the curve defined by $f_e-\vph f_o=0$ in $\IP^2$
consists of  four ordinary double points in $p_1,\ldots,p_4$.
We refer to \cite{DolgachevClassAlgGeom} and \cite{ShepherdBarronInvariantTheory} for more on the intimate relations between $\dP_5$ and genus six curves.

\subsection{The restriction map $\Pic(\dP_5)\rightarrow \Pic(C_\vph)$}\label{restrictionmap}
Let $C$ be compact Riemann surface of genus $g$, and let $\mbox{Div}(C)$ be the group of divisors on $C$. 
The Picard group of compact Riemann surface $C$ is the group of divisors on the surface modulo linear equivalence. So if $P(C)$ denotes the group of divisors of meromorphic functions, then 
$$
\Pic(C)\,=\,\mbox{Div}(C)/P(C)~.
$$
Since a divisor $D$ is a finite sum of points, with multiplicities, it has a well defined degree: 
$$
\mbox{deg}\,:\;\mbox{Div}(C)\,\longrightarrow\,\IZ\,,\qquad 
D\,=\,\sum_{p}n_pp\;\longmapsto\; \sum_p n_p~.
$$
As a meromorphic function has the same number of poles as zeroes (counted with multiplicity), one can define a subgroup $\Pic^0(C)$ of $\Pic(C)$ by:
$$
\Pic^0(C)\,:=\,\mbox{Div}^0(C)/P(C)~.
$$
By Abel's theorem, $\Pic^0(C)=Jac(C)$, the Jacobian of $C$, which is the $g$-dimensional complex torus defined as the quotient of $\IC^g$ by the period lattice, 
that is, fixing a basis $\o_1,\ldots,\o_g$ of the vector space of holomorphic 1-forms on 
$C$, the period lattice consists of the vectors $(\int_\gamma \o_1,\ldots,\int_\gamma\o_g)$ where $\gamma $ runs over all closed loops on $C$.
These groups fit together in an exact sequence:
$$
0\,\longrightarrow\,\Pic^0(C)\,\longrightarrow\,\Pic(C)\,
\stackrel{\mbox{deg}}{\longrightarrow}\, \IZ
\,\longrightarrow\,0~.
$$
As we have seen in section \ref{cover125}, a divisor $D$, whose class has order $n$ 
in $\Pic^0(C)$, so $nD$ is the divisor of a meromorphic function $f$ will define an unramified $n$:1 cover of $C$. 
As $\Pic^0(C)$ is a complex torus, it is isomorphic, as a group, 
to $(\IR/\IZ)^{2g}$. The classes $D$ with $nD=0$ are thus a subgroup isomorphic to $(\IZ/n\IZ)^{2g}$. In particular, if $C=C_\vph$ and thus $g=6$, and $n=5$ we get a subgroup $(\IZ/5\IZ)^{12}$ of five-torsion classes, whereas the subgroup of $\Pic^0(C_\vph)$ generated by the 
$D_{ij}-D_{kl}$ is a $(\IZ/5\IZ)^{3}$, 
since the covering defined by the $\sqrt[5]{g_{ij}}$ is $\tC_\vph\rightarrow C_\vph$, has degree 125 (here $g_{ij}$ has divisor $5(D_{ij}-D_{45})$ as in 
Section \ref{cover125}).

We will now identify the specific $(\IZ/5\IZ)^3\subset \Pic^0(C)$ 
which creates the covering $\tC_\vph\rightarrow C_\vph$. 
It turns that there is a quite naturally defined subgroup of $\Pic(C_\vph)$, 
which is a priori unrelated to the Dwork pencil, 
but which arises as a consequence of the special position of the curves
$C_\vph$ in $\dP_5$.

The inclusion of $C_\vph$ in the del Pezzo surface $\dP_5$ 
induces the restriction map (a homomorphism of groups)
$$
i^*:\,\Pic(\dP_5)\,\longrightarrow\,\Pic(C_\vph)~,\qquad\
\mbox{with}\quad i:\,C_\vph\,\hookrightarrow\,\dP_5~.
$$
Applying the adjunction formula, we find the canonical class on $C_\vph$:
$$
K_{C_\vph}\,=\,i^*\left(C_\vph\,+\,K_{\dP_5}\right)\,=\,
i^*\left(-K_{\dP_5}\right)
$$
where we used that the curve $C_\vph$ in $\dP_5$ has class $-2K_{\dP_5}$.
In particular, the composition $C_\vph\hookrightarrow\dP_5\hookrightarrow\IP^5$
is the canonical map. As it is an isomorphism on its image (by definition of $C_\vph$),
the curves $C_\vph$ are not hyperelliptic.

The degree two divisor $D_{ij}$ was defined as the intersection divisor of the line 
$E_{ij}\subset \dP_5$ with the curve $C_\vph\subset\dP_5$, hence
$$
D_{ij}\,=\,i^*(E_{ij})~.
$$

The group $\Pic(\dP_5)\cong\IZ^5$ has $\IZ$-basis 
$l,E_{15},\ldots,E_{45}$.
As $l=E_{ij}+E_{k5}+E_{l5}$, where $\{i,j,k,l\}=\{1,\ldots,4\}$, 
we see that the divisor $i^*l$ has degree $6$ and the 
$i^*E_{pq}=D_{pq}$ have degree two.
Thus the image of the composition of $i^*$ with 
$\mbox{deg}:\Pic(C)\rightarrow\IZ$ is the subgroup $2\IZ$ 
and the kernel of this composition is isomorphic to $\IZ^4$.
We denote this kernel by $\Pic(\dP_5)^0$:
$$
\Pic(\dP_5)^0\,:=\,
\ker(\mbox{deg}\circ i^*)\,=\oplus_{i=1}^4\IZ\alpha_i~,
$$
where the $\IZ$-basis $\alpha_i$ of $\Pic(\dP_5)^0$ is defined by
$$
\alpha_1\,=\,E_{15}-E_{25},\quad \alpha_2\,=\,E_{25}-E_{35},\quad
\alpha_3\,=\,E_{35}-E_{45},\quad \alpha_4\,=\,l-E_{15}-E_{25}-E_{35}~.
$$
As we have  $C_\vph= -2K_{\dP_5}$ in $\Pic(\dP_5)$, the divisors on $\dP_5$ which intersect $C_\vph$ in a divisor of degree $0$ form the subgroup  
$K_{\dP_5}^\perp$, so 
$$
\Pic(\dP_5)^0\,=\,K_{\dP_5}^\perp\,=\,
\{\,D\,\in\,\Pic(\dP_5)\,:\,D\cdot K_{\dP_5}\,=\,0\,\}~.
$$
We recall the well-known fact that the intersection matrix of the $\alpha_i$ is the Cartan matrix of the root system $A_4$, up to sign:
$$
(\alpha_i,\alpha_j)\,=\,\left\{\begin{array}{rcl} -2&\mbox{if}&i=j~,\\
1&\mbox{if}&|i-j|=1~,\\ 0&\mbox{else}~.&\end{array}\right.
$$

We now have a homomorphism, induced by $i^*$, 
but denoted by the same symbol,
$$
i^*:\,\Pic(\dP_5)^0\,\longrightarrow \,Pic^0(C)~.
$$
As the $\alpha_i$, $i=1,2,3$, are of the form $E_{ij}-E_{pq}$, 
their images $i^*(\alpha_i)$ are of the form $D_{ij}-D_{pq}$
which are elements of order five in $\Pic^0(C)$.
Finally, using that $l=E_{34}+E_{15}+E_{25}$, we see that
$$
\alpha_4\,=\,l\,-\,E_{15}\,-\,E_{25}\,-\,E_{35}\,=\,
E_{34}\,+\,E_{15}\,+\,E_{25}\,-\,(E_{15}\,+\,E_{25}\,+\,E_{35})\,=\,
E_{34}\,-\,E_{35}~,
$$
hence $i^*(\alpha_4)=D_{34}-D_{45}$ is also 5-torsion. 
The image of $i^*$ is generated by the classes $i^*(\alpha_j)$:
{\renewcommand{\arraystretch}{1.5}
$$
\begin{array}{rcl}
\mbox{im}(i^*)&=&
\langle D_{15}\,-\,D_{25},\, D_{25}\,-\,D_{35},\,D_{35}-D_{45},\,
D_{34}\,-\,D_{35}\,\rangle\\
&=&\langle D_{15}\,-\,D_{45},\, D_{25}\,-\,D_{45},\,D_{35}-D_{45},\,
D_{34}\,-\,D_{45}\,\rangle~.
\end{array}
$$
}
Thus $\mbox{im}(i^*)\cong(\IZ/5\IZ)^n$ for some $n\leq 4$.
There is a further relation between these classes, given by the divisor of the function $k_{14}/l_1l_2$ 
(notice that we take the quotient of two polynomials of bidegree $(2,2)$, so this quotient is a well-defined meromorphic function on $C^0_\vph$). We use Table \ref{Divs}
to find the divisor $(k_{14})$ of $k_{14}$ in $\mbox{Div}(C_\vph)$ (and we  simply write $l$ for $i^*(l)$):
{\renewcommand{\arraystretch}{1.5}
$$
\begin{array}{rcl}
(k_{14})&:=&
D_{23}+D_{25}+3D_{24}+3D_{12}\\
&=&(l-D_{15}-D_{45})+D_{25}+3(l-D_{15}-D_{35})+3(l-D_{35}-D_{45})\\
&=&7l-4D_{15}+D_{25}-6D_{35}-4D_{45}\\
&=&7(D_{34}+D_{15}+D_{25})-4D_{15}+D_{25}-6D_{35}-4D_{45}\\
&=&7D_{34}+3D_{15}+8D_{25}-6D_{35}-4D_{45}~.
\end{array}
$$
}
Similarly, the divisor of $l_1$ is:
{\renewcommand{\arraystretch}{1.5}
$$
\begin{array}{rcl}
(l_1)&=&D_{13}+D_{23}+D_{34}+D_{35}\\
&=&(l-D_{25}-D_{45})+(l-D_{15}-D_{45})+D_{34}+D_{35}\\
&=&2l-D_{15}-D_{25}+D_{35}-2D_{45}+D_{34}\\
&=&2(D_{34}+D_{15}+D_{25})-D_{15}-D_{25}+D_{35}-2D_{45}+D_{34}\\
&=&3D_{34}+D_{15}+D_{25}+D_{35}-2D_{45}~.
\end{array}
$$
}
As $(l_2)=D_{15}+D_{25}+D_{35}+D_{45}$, the linear equivalence $(k_{14})=(l_1)+(l_2)$ gives the following relation
in $\Pic(C_\vph)$:
$$
7D_{34}+3D_{15}+8D_{25}-6D_{35}-4D_{45}
\,=\,
3D_{34}+2D_{15}+2D_{25}+2D_{35}-D_{45}~.
$$
Using $5D_{25}=5D_{35}$ in $\Pic(C_\vph)$ we get:
$$
4D_{34}\,=\,-D_{15}-6D_{25}+8D_{35}+3D_{45}\,=\,
-D_{15}-D_{25}+3D_{35}+3D_{45}~\qquad(\in \Pic(C_\vph))~.
$$
Now we write $4D_{34}=-D_{34}+5D_{45}$ and use $5D_{34}=5D_{45}$ in $\Pic(C_\vph)$ to obtain 
$$
D_{34}\,=\,D_{15}+D_{25}-3D_{35}+2D_{45}~\qquad(\in \Pic(C_\vph))~.
$$
This gives the following relation in $\Pic^0(C_\vph)$:
$$
D_{34}-D_{45}\,=\,(D_{15}-D_{45})\,+\,(D_{25}-D_{45})\,+\,2(D_{35}-D_{45})~.
$$
Therefore $\mbox{im}(i^*)\subset \Pic^0(C_\vph)$ 
can be generated by three elements and thus $n\leq 3$.

We give a table which gives the$(a_1,a_2,a_3)\in(\IZ/5\IZ)^3$ such that 
the classes $e_{ij}:=D_{ij}-D_{45}=i^*(E_{ij}-E_{45})$ 
can be written as $a_1e_{15}+a_2e_{25}+a_3e_{35}$.
$$
\begin{array}{rclrclrcl}
e_{15}&\mapsto&(1,0,0),   & e_{25}&\mapsto&(0,1,0), &
e_{35}&\mapsto&(0,0,1),\\
e_{12}&\mapsto&(2,2,1),  &
e_{13}&\mapsto&(2,1,2), &e_{23}&\mapsto&(1,2,2),  \\ 
e_{14}&\mapsto&(2,1,1),&\quad e_{24}&\mapsto&(1,2,1),&\quad e_{34}&\mapsto&(1,1,2)~.\\
\end{array}
$$

As we have seen, $5D_{ij}=5D_{pq}$ for any indices. 
Thus we have a rather peculiar divisor class of degree $5\cdot 2=10$ in $\Pic(C_\vph)$. This is actually the canonical class on $C_\vph$. In fact,
{\renewcommand{\arraystretch}{1.5}
$$
\begin{array}{rcl}
K_{C_\vph}&=&i^*(-K_{\dP_5}) \\
&=&i^*(3l-E_{15}-E_{25}-E_{35}-E_{45})\\
&=& 3(D_{34}+D_{15}+D_{25})-D_{15}-D_{25}-D_{35}-D_{45}\\
&=&3(2D_{15}+2D_{25}-3D_{35}+2D_{45})-D_{15}-D_{25}-D_{35}-D_{45}\\
&=&5D_{15}+5D_{25}-10D_{35}+5D_{45}\\
&=&5D_{15}
\end{array}
$$
}
where in the last step we used that $5D_{ij}=5D_{15}$.

To show that $n\geq 3$ and thus $n=3$, we use the $\cA_5$-action on the
subgroup $\mbox{im}(i^*)\cong(\IZ/5\IZ)^n$ of  $\Pic^0(C_\vph)$. 
As $\cA_5$ is a simple group, the image of $\cA_5$ in $\Aut(\mbox{im}(i^*))$
is either the identity or isomorphic to $\cA_5$. 
In the first case, by applying $(23)(45)\in \cA_5$ to $D_{12}-D_{45}$ 
we would get that $D_{12}-D_{45}=D_{13}-D_{45}$ in $\Pic^0(C_\vph)$ 
and hence that $D_{12}-D_{45}$ is the divisor of a meromorphic function, which is not the case as $C_\vph$ is not hyperelliptic.
Thus we obtain an injective homomorphism $\cA_5\rightarrow GL(n,\IZ/5\IZ)$.
If $n=1$, this is impossible as $\sharp \cA_5 > \sharp GL(1,\IZ/5\IZ)=4$.
If $n=2$, we consider the action of the subgroup 
$\{e,(12)(34),(13)(24),(14)(23)\}$, which is isomorphic to $(\IZ/2\IZ)^2$. 
The eigenvalues of $(12)(34)$ on $(\IZ/5\IZ)^2$ are $1,1$; $1,-1$ or $-1,-1$.
The first case is impossible since the homomorphism is injective and an automorphism with eigenvalues $1,1$ which is not the identity has order $5$.
The last case is also impossible since then either $(12)(34)$ should commute with all other elements of $\cA_5$ or the element would have order five again.
Thus the eigenvalues are $1,-1$ and we can diagonalize the automorphism. The same is true for the other two non-trivial elements in the subgroup. Since the subgroup is commutative, these three automorphisms can be diagonalized on the same basis. But then one of the automorphisms must be $-I$, which commutes with any other automorphism, again a contradiction. 
Therefore $n\geq 3$. 

We conclude that 
$\mbox{im}(i^*)\subset \Pic^0(C_\vph)$ consists of the divisor classes which create the unramified 125:1 covering $\tC_\vph\rightarrow C_\vph$. 
So for any non-zero
$D\in \mbox{im}(i^*)$ there is a meromorphic function $f_D$ on $C_\vph$ 
with divisor $5D$ and the Riemann surface $C_D$ defined by the polynomial 
$T^5-f_D$ is a 5:1 unramified cover of $C_\vph$ which fits in a diagram
$\tC_\vph\rightarrow C_D\rightarrow C_\vph$. The fiber product over $C_\vph$ of three suitably choosen $C_D$ will be isomorphic to $\tC_\vph$.
Since the image of the curve $C_\vph$ in the canonical embedding lies in a unique del Pezzo surface of degree 5 (see \cite{ShepherdBarronInvariantTheory}), we have the remarkable fact that the curves $\tC_\vph$ which parametrize the lines in Dwork pencil are
intrinsically determined by the curves $C_\vph$ in the Wiman pencil.

In this section we have shown that $i^*$ maps $\Pic(\dP_5)^0\cong Q(A_5)$
onto $(\IZ/5\IZ)^3$. In terms of root systems, this map is well-known. 
The root lattice $Q(A_5)$ is a sublattice, of index $5$, of the weight lattice $P(A_5)$. Thus $5P(A_5)$ is a sublattice, of index $125$, of $Q(A_5)$ 
and thus $Q(A_5)/5P(A_5)\cong (\IZ/5\IZ)^3$. 
The map $i^*$ can be identified with the quotient map 
$Q(A_5)\rightarrow Q(A_5)/5P(A_5)$.
\vskip2in
{\bf\large Acknowledgements}
\vskip10pt
We wish to thank S. Katz for discussions. PC and XD also wish to thank the Perimeter Institute and ICTP, Trieste, for support and hospitality while they were engaged on this this~project.
\newpage
\raggedright
\bibliographystyle{utphys}
\bibliography{BibliographyLines}

\providecommand{\href}[2]{#2}\begingroup\raggedright\begin{thebibliography}{10}

\bibitem{Schubert1}
H.~Schubert, ``Das {C}orrespondenzprincip f\"{u}r {G}ruppen von $n$ {P}unkten
  und von $n$ {S}trahlen,'' {\em Math. Annalen} {\bfseries 12} (1877) 180--201.

\bibitem{Schubert2}
H.~Schubert, ``Die $n$-dimensionale {V}erallgemeinerung der {A}nzahlen f\"{u}r
  die vielpunktig ber\"{u}hrenden {T}angenten einer punktallgemeinen
  {F}l\"{a}che $m$-ten {G}rades,'' {\em Math. Annalen} {\bfseries 26} (1886)
  52--73.

\bibitem{CayleyLines}
A.~Cayley, ``On the triple tangent planes of surfaces of the third order,''
  {\em Camb. and Dublin Math. Journal} {\bfseries IV} (1849) 118--132.

\bibitem{SalmonLines}
G.~Salmon, ``On the triple tangent planes of surfaces of the third order,''
  {\em Camb. and Dublin Math. Journal} {\bfseries IV} (1849) 152--260.

\bibitem{HendersonLines}
A.~Henderson, {\em The twenty-seven lines upon the cubic surface}.
\newblock No.~13 in {C}amb. {T}racts. Cambridge University Press, 1911.

\bibitem{Jordan1870}
C.~Jordan, {\em Trait{\'e} des substitutions et des {\'e}quations
  alg{\'e}braiques}.
\newblock Gauthier-Villars, Paris, 1870.

\bibitem{HarrisGaloisGroups}
J.~Harris, ``Galois groups of enumerative problems,'' {\em Duke {M}ath. {J}.}
  {\bfseries 46} (1979) 685--724.

\bibitem{MR1024767}
A.~Albano and S.~Katz, ``Lines on the {F}ermat quintic threefold and the
  infinitesimal generalized {H}odge conjecture,'' {\em Trans. Amer. Math. Soc.}
  {\bfseries 324} no.~1, (1991) 353--368.

\bibitem{KatzDegenerations}
S.~Katz, ``Degenerations of quintic threefolds and their lines,'' {\em Duke
  {M}ath. {J}.} {\bfseries 50} no.~4, (1983) 1127--1135.

\bibitem{MR1085631}
A.~Albano and S.~Katz, ``van {G}eemen's families of lines on special quintic
  threefolds,'' {\em Manuscripta Math.} {\bfseries 70} no.~2, (1991) 183--188.

\bibitem{Mustata:fk}
A.~{\Mustata}, ``Degree 1 curves in the {D}work pencil and the mirror
  quintic,'' \href{http://dx.doi.org/10.1007/s00208-011-0668-x}{{\em Math.
  Annalen} (2010) 1--34},
  \href{http://arxiv.org/abs/math/0311252v1[math.AG]}{{\ttfamily
  arXiv:math/0311252v1[math.AG]}}.

\bibitem{Wiman}
A.~Wiman, ``{Z}ur {T}heorie der endlichen {G}ruppen von birationalen
  {T}ransformationen in der {E}bene,'' {\em Math. Annalen} {\bfseries 48}
  (1897) 195--240.

\bibitem{Edge}
W.~L. Edge, ``A pencil of four-nodal plane sextics,'' {\em Math. Proc. Camb.
  Phil. Soc.} {\bfseries 89} (1981) 413--421.

\bibitem{Candelas:2004sk}
P.~Candelas, X.~de~la Ossa, and F.~Rodriguez~Villegas, ``{Calabi-Yau manifolds
  over finite fields II},''
\href{http://arxiv.org/abs/hep-th/0402133}{{\ttfamily arXiv:hep-th/0402133
  [hep-th]}}.

\bibitem{ForsterRiemannSurfaces}
O.~Forster, {\em Lectures on Riemann surfaces}.
\newblock No.~81 in Graduate Texts in Mathematics. Springer-Verlag, New York,
  1991.

\bibitem{SchoenMumfordHorrocksBundle}
C.~Schoen, ``On the geometry of a special determinantal hypersurface associated
  to the {M}umford-{H}orrocks vector bundle,'' {\em J. Reine Angew. Math.}
  {\bfseries 364} (1986) 85--111.

\bibitem{HigherTranscendentalFunctions}
A.~Erd{\'e}lyi, W.~Magnus, F.~Oberhettinger, and F.~G. Tricomi, {\em Higher
  Transcendental Functions}.
\newblock Mc{G}raw-{H}ill {B}ook {C}ompany, 1955.

\bibitem{Forsyth}
A.~R. Forsyth, {\em Theory of functions of a complex variable}.
\newblock Cambridge University Press, 1893.

\bibitem{DolgachevClassAlgGeom}
I.~V. Dolgachev, {\em Classical Algebraic Geometry: a modern view}.
\newblock Cambridge University Press, to appear.

\bibitem{ShepherdBarronInvariantTheory}
N.~Shepherd-Barron, ``Invariant theory for ${S}_5$ and the rationality of
  ${M}_6$,'' {\em Compositio Math.} {\bfseries 70} (1989) 13--25.

\end{thebibliography}\endgroup
\end{document}